\documentclass{paper}
\usepackage{amssymb}
\usepackage{amsmath}
\usepackage{amsfonts}
\usepackage{graphicx}
\usepackage{wrapfig}
\usepackage{enumerate}
\usepackage{mathrsfs}
\usepackage{hyperref}
\usepackage{algorithm}
\usepackage{algpseudocode}
\usepackage{fancyhdr}
\usepackage{color}
\usepackage{amsthm}
\newtheorem{remark}{Remark}
\renewcommand{\vec}[1]{\ensuremath{\boldsymbol{#1}}}

\newcommand{\Reals}{\ensuremath{{\mathbb{R}}}}

\newcommand{\diff}{\ensuremath{{\operatorname{d}}}}
\renewcommand{\O}{\ensuremath{{\Omega}}}

\newcommand{\normal}{\ensuremath{{\vec{\nu}}}}

\newcommand{\G}{\ensuremath{{\Gamma}}}

\newcommand{\Lp}[1]{\ensuremath{\operatorname{L}_{#1}}}

\usepackage{upgreek}

\newcommand{\p}{\ensuremath{{\mathsf{p}}}}

\usepackage{color}

\definecolor{MyGreen}{rgb} {0.05,0.4,0.05}
\definecolor{RedViolet}{rgb} {0.1,0.1,0.75}

    \providecommand{\Ignore}[1]{}
    \providecommand{\ignore}[1]{}
    \providecommand{\freeze}[1]{}
    \providecommand{\crossout}[1]{{\textcolor{red!20}{#1}}}
    \providecommand{\highlight}[1]{{\color{red}#1}}
    \providecommand{\standout}[1]{\colorbox{red}{\textcolor{green}{#1}}}
    
    \provideboolean{shownotes}
    \setboolean{shownotes}{true}
    \newcounter{margnote}[page]
    \providecommand{\margnotemark}{{\standout{\upshape\texttt{\arabic{margnote}}}}}
    \providecommand{\margnote}[2][]{
      \ifthenelse{
        \boolean{shownotes}
      }{
        \stepcounter{margnote}
        \margnotemark
        \marginpar{
          \texttt{
 \begin{minipage}{2.5cm}
              \raggedright\tiny
              \margnotemark
              {\ifx& #1&{}\else{#1 says:}\fi} 
              #2
            \end{minipage}
          }
        }
      }{
      }
    }
    \providecommand{\mathnote}[2][]{
      \ifthenelse{
        \boolean{shownotes}
      }{
        \stepcounter{margnote}
        \text{
          \standout{
            \texttt{
              \tiny
              \margnotemark
              #1:
              #2
            }
          }
        }
      }{
      }
    }
    \provideboolean{showtodo}

    \providecommand{\todo}[1]{\ifthenelse{\boolean{showtodo}}{\margnote[To do.]{#1}}{}}
    \providecommand{\Todo}[1]{
      \ifthenelse{\boolean{showtodo}}{
        \begin{center}
        \begin{tikzpicture}
         \node[fill=a!17]{
           \begin{minipage}{\textwidth}
             \texttt{\bfseries{\small #1}}
           \end{minipage}
         };
        \end{tikzpicture}
        \end{center}
      }{}}

    \provideboolean{showcomments}
    \providecommand{\margincomment}[1]{
    \ifthenelse{\boolean{showcomments}}{\marginpar{\tiny #1}}{}
    }
    \provideboolean{showchanges}
    \setboolean{showchanges}{true}
    \providecommand{\changes}[1]{
      \ifthenelse{\boolean{showchanges}}{{\highlight{#1}}}{#1}
    }
    \providecommand{\changefromto}[3][replace with]{
      \ifthenelse{\boolean{showchanges}}
      {{\crossout{#2}\margnote{#1}}{\highlight{#3}}}
      {#3\xspace}
    }
    \providecommand{\ChangePar}[2]{
      \ifthenelse{\boolean{showchanges}}
      {{\par$\mapsfrom$ \textcolor{red!20}{#1}}{\par$\mapsto$ \textcolor{blue}{#2}}}
      {\par #2}
    }
    \providecommand{\InsertPar}[1]{
      \ifthenelse{\boolean{showchanges}}
      {{\par$\mapsto$ \textcolor{blue}{#1}}}
      {\par #1}
    }

\usepackage[toc,page]{appendix}
\numberwithin{equation}{section}
\begin{document}

\pagestyle{fancy}
\lhead{optimal control paper, Section \thesection}
\rhead{\thepage}


\newcommand{\Eqref}[1]{(\ref{#1})}




\title{{A robust and efficient adaptive multigrid solver for the optimal control of phase field formulations of geometric evolution laws}}



\author{F. Yang}{F.W.Yang@sussex.ac.uk}{a}

\author{C. Venkataraman}{C.Venkataraman@sussex.ac.uk}{a}

\author{V. Styles}{V.Styles@sussex.ac.uk}{a}

\author{A. Madzvamuse}{A.Madzvamuse@sussex.ac.uk}{a}

\affiliation{a}{Department of Mathematics}
  {University of Sussex}


\begin{abstract}
We propose and investigate a novel solution strategy to efficiently and accurately compute approximate solutions to semilinear optimal control problems, focusing on the optimal control of phase field formulations of geometric evolution laws.
The optimal control of geometric evolution laws arises in a number of applications in fields including material science, image processing, tumour growth and cell motility.
Despite this, many open problems remain in the analysis and approximation of such problems.
In the current work we focus on a phase field formulation of the optimal control problem, hence exploiting the well developed mathematical theory for the optimal control of semilinear parabolic partial differential equations.
Approximation of the resulting optimal control problem is computationally challenging, requiring massive amounts of computational time and memory storage.
The main focus of this work is to propose, derive, implement and test an efficient solution method for such problems.
The solver for the discretised partial differential equations is based upon a geometric multigrid method incorporating advanced techniques to deal with the nonlinearities in the problem and utilising adaptive mesh refinement.
An in-house two-grid solution strategy for the forward and adjoint problems, that significantly reduces memory requirements and CPU time, is proposed and investigated computationally.
Furthermore, parallelisation as well as an adaptive-step gradient update for the control are employed to further improve efficiency.
Along with a detailed description of our proposed solution method together with its implementation we present a number of computational results that demonstrate and evaluate our algorithms with respect to accuracy and efficiency.
A highlight of the present work is simulation results on the optimal control of phase field formulations of geometric evolution laws in 3-D which would be computationally infeasible without the solution strategies proposed in the present work.
 
\end{abstract}

%
%

\section{Introduction}
\label{sec_intro}

The optimal control of geometric evolution equations or more generally free boundary problems arises in a number of applications.
In image processing the tracking of deformable objects may be formulated as the optimal control of a suitably chosen evolution law  \cite{papadakis2007variational}. 
A number of applications arise from problems in material science such as the control of nanosturcture  through electric fields \cite{Hauber1,Hauber2}. An important and topical application area
is the image driven modelling of biological processes, such as tumour growth \cite{hogea2008image} or cell migration \cite{2013arXiv1311.7602C}, in which parameters (or functions) in a model are estimated from experimental imaging data. 
In a recent study we proposed an optimal control approach to whole cell tracking \cite{Costas}, i.e., the reconstruction of whole cell morphologies in time from a set of static images, in which the cell tracking problem was formulated as the optimal control of a geometric evolution equation \cite{Costas}. In general the approximation of such optimal control problems is computationally intensive both in terms of central processing unit (CPU) time and memory. Hence, the development of robust and efficient solvers for such problems with a view to reducing CPU time (or simply wall-clock time) and memory requirements is a worthwhile research direction.

In the current work we consider the optimal control of geometric evolution laws of forced mean curvature flow type. We denote by $\G(t)$,
 a closed oriented smoothly evolving $d-1$ dimensional hypersurface in $\Reals^d, d=2,3$ with outward pointing unit normal $\normal$. The motion of $\G(t)$ satisfies a volume constrained mean curvature flow with forcing, i.e., given an initial surface $\G(0)$, the velocity $\vec V$ of $\G$ is given by
 \begin{equation}\label{eqn:vc_mcf}
 \vec V(\vec x,t)=\left(-\sigma H(\vec x,t)+\eta(\vec x,t)+\lambda_V(t)\right)\normal(\vec x,t)\quad \vec x\in\G(t), t\in(0,T],
 \end{equation}
 where $\sigma>0$ represents the surface tension, $H$ denotes the mean curvature (which we take to be the sum of the principal curvatures) of $\G$, $\eta$ is a space time distributed forcing and $\lambda_V$ is a spatially uniform Lagrange multiplier enforcing volume constraint. We assume we are given an initial interface $\G^0$ and a target interface $\G_{obs}$ both of which are smooth closed oriented $d-1$ dimensional hypersurfaces.
 
The optimal control problem, which is the focus of the current work, consists of finding a space time distributed forcing $\eta$ in \eqref{eqn:vc_mcf} such that with $\G(0)=\G^0$, the interface position at time $T$ corresponding to the solution of \eqref{eqn:vc_mcf}, $\G(T)$, is ``close'' to the observed data $\G_{obs}$. 
We have deliberately refrained from stating precisely what is meant by $\G(T)$ being close to $\G_{obs}$ as in the sharp interface setting it is not obvious what constitutes a good choice of metric to measure the difference between two surfaces.
In particular standard measures such as the Haussdorff distance are typically non-smooth and this complicates the approximation of the optimal control problem.
Moreover, the theory of optimal control of geometric evolution laws is in its infancy, in fact only recently has progress been made on the optimal control of parabolic equations on evolving surfaces even in the case of prescribed evolution \cite{vierling2014parabolic}.
On the other hand, the theory for the optimal control of semilinear parabolic equations is more mature (see, for example, \cite{Troltzsch}) and if one considers a diffuse interface representation of the surfaces then standard measures of distance may be used (e.g., $\Lp{2}$).
In light of the above, we consider the phase field approximation of (\ref{eqn:vc_mcf}) given by the volume constrained Allen-Cahn equation, see \eqref{equ_Allen-Cahn}.
We approximate our initial and target data $\G^0$ and $\G_{obs}$ by diffuse interface representations $\phi_0$ and $\phi_{obs}$ respectively (see for example \cite{2013arXiv1311.7602C} for details on how to construct such representations). 

Our strategy for approximating the solution to the optimal control problem consists of  an iterative adjoint based solution method, c.f., Sections \ref{sec_oc} and \ref{sec_numerical}. The method is particularly computationally intensive for a number of reasons.
\begin{enumerate}
\item The iterative adjoint based solution method we employ necessitates multiple solves of the forward and adjoint problems.
\item As the state equation is of Allen-Cahn type grid adaptivity for the solution of the forward equation is mandatory (at least in 3-D). 
\item The computed state enters the adjoint equation which is posed backwards in time. Hence the state equation must first be solved over the whole time interval with the computed states  stored and then the adjoint equation is solved backwards in time. Thus the algorithm requires large amounts of data storage.
\item We want to consider small values of the interfacial width parameter as many of the applications from cell biology that we are interested in involve interfaces with large curvatures and small scale features  which we wish to resolve with our diffuse interface approximation. This imposes strong restrictions on the grid for the solution of the state problem.
\end{enumerate}
In this work we  focus on developing a robust and efficient solver for the problem. We employ a fast parallel adaptive multigrid solution method for the forward equation. The use of adaptive grids and parallelisation allows us to compute with relatively small values of the interfacial width parameter. For the adjoint equation we make the observation that as the PDE is linear it may be possible to relax the restrictions on the grid needed for the solution of the state equation. Hence we employ a parallel multigrid solver for the adjoint equation on a uniform grid that is typically coarser than the adaptive grid used for the approximation of the state. We also consider a simple adaptive strategy for our iterative steepest descent based algorithm for the update of the control.

Major findings in the present work include that the two-grid solution strategy we propose in which the state equation is solved on an adaptive grid and the adjoint problem is solved on a coarse uniform grid appears to have only a minor detrimental impact 
on accuracy whilst the savings in terms of memory and CPU time are considerable as the state is only stored on the coarse grid.
We also propose an adaptive algorithm for the iterative update which finds the optimal control $\eta$.
We expect that our findings are of relevance beyond the optimal control of Allen-Cahn equations alone and this could include the development of efficient schemes for the optimal control of semilinear parabolic equations in general.

The remainder of this article is structured as follows.
In Section \ref{sec_oc} we present the optimal control problem and formally derive the optimality conditions used in the algorithm for its approximation. In Section \ref{sec_numerical} we
describe our solver for the optimal control problem, the key import of this work. Firstly, in Section \ref{sec_31} we summarise the procedures required for solving this optimal control problem. We also outline the fully discrete methods for the approximation of the state and adjoint equations in Section \ref{sec_32}. 
We discuss the two-grid solution strategy in Section \ref{sec_33}; and the adaptive-$\alpha$ algorithm in Section \ref{adaptive-a} which improves the control update.
Section \ref{sec_results} contains results of our numerical experiments with the implemented solver. 
We use a 2-D benchmark problem to demonstrate the convergence of the proposed model, our multigrid performance where a linear complexity is shown, effectiveness of the adaptive-$\alpha$ algorithm and two-grid strategy.
We use a 3-D benchmark problem to illustrate the importance of the proposed two-grid solution strategy in terms of the saved memory spaces.
We also show 2-D and 3-D irregular shapes.
Finally in Section \ref{sec_concolusions} we summarise our major findings and discuss directions for future work.

\section{Optimal control of a forced Allen-Cahn equation with volume constraint}
\label{sec_oc}

As outlined in Section \ref{sec_intro} the prototype state equation we consider in this work consists of the volume constrained Allen-Cahn equation with forcing, 
\begin{equation}
\label{equ_Allen-Cahn}
\begin{cases}
\epsilon \frac{\partial}{\partial t} \phi(\vec x,t) = \epsilon \triangle\phi(\vec x,t) - \epsilon^{-1}G'(\phi(\vec x,t))+ \eta(\vec x,t) + \lambda(t)\quad &\text{in }\O\times(0,T],\\
\phi(\vec x,0)=\phi^0(\vec x)\quad&\text{in }\O,\\
\nabla\phi(\vec x,t)\cdot\normal_\O(\vec x)=0\quad&\text{on }\partial\O,
\end{cases}
\end{equation}
where $\phi(\vec x,t)$ is the phase field variable, $\epsilon > 0$ is the parameter governing the interfacial width of the diffuse interface, $G(\eta) = \frac{1}{4}(1-\eta^2)^2$ is a double well potential which has minima at $\pm1$ and $\lambda$ is a time-dependent constraint on the mass that models a volume constraint \cite{Blowey} and $\normal_\O$ is the normal to $\partial\O$.

As the volume enclosed by the target and initial interfaces may differ, i.e., $\int_\O\phi^0\neq\int_\O\phi_{obs}$, enforcing conservation of mass is inappropriate instead we proceed as in \cite{Costas} and enforce a constraint on the linear interpolant of the mass of the initial and target diffusive interfaces. To this end we define
\begin{equation}
\label{equ_M}
M_{\phi}(t) := \int_{\Omega}\left[\phi^0(\vec x) + \frac{t}{T}\left(\phi_{obs}(\vec x) - \phi^0(\vec x)\right) \right],
\end{equation}
and the volume constraint $\lambda(t)$ in (\ref{equ_Allen-Cahn}) is then determined such that for $t\in(0,T]$
\begin{equation}
\int_{\Omega}\phi(\vec x,t)=  M_\phi(t).
\end{equation}

In order to formulate our optimal control problem we introduce the objective functional $J$, which we seek to minimise
\begin{equation}
\label{equ_J}
J(\phi,\eta) = \frac{1}{2}\int_{\Omega}\left(\phi(\vec x,T) - \phi_{obs}(\vec x) \right)^2\diff \vec x + \frac{\theta}{2} \int^{T}_{0}\int_{\Omega}\eta^2(\vec x,t)\diff \vec x\diff t,
\end{equation}
where $\theta > 0$ is a regularisation parameter.
The first term of the right-hand side of (\ref{equ_J}) is the so called fidelity term which measures the distance between the solution of the model and the target data $\phi_{obs}$ and the second term is the regularisation which is necessary to ensure a well-posed problem \cite{Troltzsch}.

The optimal control problem we consider in this work may now be stated as the following minimisation problem.  Given initial data $\phi^0$ and target data $\phi_{obs}$,  find a space-time distributed forcing $\eta^*:\O\times[0,T)\to\Reals$ such that with $\phi$ a solution of \eqref{equ_Allen-Cahn} with initial condition $\phi(\cdot,0)=\phi^0(\cdot)$, the forcing $\eta^*$ solves the minimisation problem
\begin{equation}
\label{equ_minimisation}
\text{min}_{\eta} J(\phi,\eta), \;\; \text{with $J$ given by (\ref{equ_J})}.
\end{equation}

In order to apply the theory of optimal control of semilinear PDEs for the solution of the minimisation problem, we briefly outline the derivation of the optimality conditions, for further details see for example \cite{Troltzsch,hinze2009optimization}. 
Adopting a Lagrangian approach, we introduce  the Lagrange multiplier  (adjoint state) $p$ and the Lagrangian functional $\mathcal{L}(\phi,\eta,p)$ defined by
\begin{equation}
\label{equ_Lagrangian}
\mathcal{L}(\phi,\eta,p) = J(\phi,\eta) - \int^T_0 \int_\Omega \left(\epsilon \frac{\partial}{\partial t} \phi - \epsilon \triangle\phi+ \epsilon^{-1}G'(\phi)+ \eta + \lambda \right) p.
\end{equation}

Assuming the existence of an optimal control $\eta^*$ and the associated optimal state $\phi^*$ and requiring stationarity of the Lagrangian at $(\phi^*,\eta^*)$  yields the (formal) first order optimality conditions \cite{Troltzsch,hinze2009optimization}
\begin{eqnarray}
\delta_\phi \mathcal{L}(\phi^*,\eta^*,p)\phi = 0, &\;\;\forall \phi : \phi(\vec x,t=0) = 0, \label{equ_for_p}\\
\delta_\eta \mathcal{L}(\phi^*,\eta^*,p)\eta = 0, &\;\;\forall \eta.  \label{equ_for_eta}
\end{eqnarray}
Condition (\ref{equ_for_p}) yields the linear parabolic adjoint equation posed backwards in time
\begin{equation}
\label{equ_adjoint}
\begin{cases}
\frac{\partial}{\partial t}p(\vec x,t) = - \triangle p(\vec x,t) + \epsilon^{-2}{G''\left(\phi\left(\vec x,t\right)\right)}p(\vec x,t)\quad&\text{in }\O\times[0,T),\\
p(\vec x,T)=\phi(\vec x,T)-\phi_{obs}(\vec x)\quad&\text{in }\O,\\
\nabla p(\vec x,t)\cdot\normal_\O(\vec x)=0\quad&\text{on }\partial\O\times[0,T).
\end{cases}
\end{equation}
Condition (\ref{equ_for_eta}) together with the Riesz representation theorem yields the optimality condition for the control \cite{Troltzsch}
\begin{equation}\label{eqn:eta_opt}
\delta_\eta \mathcal{L}(\phi^*,\eta^*,p) = \theta \eta^* + \frac{1}{\epsilon}p = 0.
\end{equation}

We note that our approach to the optimal control problem involving the formulation of the adjoint problem appears to require a smooth potential $G$. The formulation of the adjoint problem is to our best knowledge an open problem for other widely used, but non smooth or unbounded, potentials such as the obstacle or logarithmic potential.

\section{Numerical solution methods}
\label{sec_numerical}

In this section, we outline the numerical methods for solving the proposed optimal control of geometric evolution laws.
For ease of exposition, we separate this complex solution procedure into three parts.
The first part, which deals with the update of the optimal control $\eta$, is discussed in Section \ref{sec_31}, with the assumption that the solutions of $\phi$ and $p$ have already been obtained.
In Section \ref{sec_32}, we describe the second part that involves the spatial and temporal discretization schemes for the forward ($\phi$) and adjoint ($p$) problems, i.e., the phase field Allen-Cahn equation (\ref{equ_Allen-Cahn}) with volume constraint and the adjoint equation (\ref{equ_adjoint}) respectively.
In order to solve this problem efficiently and accurately, we employ several state-of-the-art algorithms, described in Section \ref{sec_33}, including an important in-house two-grid solution strategy, parallelised multigrid solution methods as well as dynamic adaptive mesh refinement.
This forms the third and the last part of our solution procedure.

\subsection{Adaptive iterative update for the optimal control}
\label{sec_31}\label{adaptive-a} 
The control $\eta$ is updated and obtained through an iterative approach, where the Allen-Cahn and adjoint equations (\ref{equ_Allen-Cahn}) and (\ref{equ_adjoint}) respectively,  for each fixed time frame $[0,T]$ have to be solved repeatedly.
The computational requirement for finding a satisfactory $\eta$ is large.
This is mainly due to two reasons; first the necessity of repeatedly solving both Allen-Cahn and adjoint equations sequentially and second, the state ($\phi$) and the forcing ($\eta$) must be stored at all iterations.

We denote a superscript $\ell$ for the $\eta$ iteration, and at $\ell=0$, we take
\begin{equation}
\eta^{\ell=0} = 0 \;\; \text{on}\;\; \Omega \times [0, T)
\end{equation}
as our initial guess for the control.
A better initial guess for the control may be necessary for certain examples or applications, however for the benchmark examples presented in this paper the simple constant zero initial guess stated above was sufficient for convergence.

For the purpose of demonstration, let us assume that both the state and adjoint equations are solved by some known method with an acceptable accuracy.
With this assumption, a gradient-based iterative update of the control, following the steepest descent approach, is employed using \eqref{eqn:eta_opt} and the update is given by
\begin{equation}
\label{equ_steepest_descent}
\eta^{\ell+1} = \eta^{\ell} - \alpha\left(\theta \eta^{\ell} + \frac{1}{\epsilon}p^{\ell} \right), \;\; \text{on} \;\; \Omega \times [0,T),
\end{equation}
where $\ell+1$ denotes the next $\eta$ iteration and $\ell$ indicates the current $\eta$ iteration.

The whole procedure is repeated until the objective function $J$ (see \eqref{equ_J}) satisfies some pre-defined tolerances.
There are two criteria: an absolute criterion and a relative criterion.
The former terminates when the obtained $J$ is smaller than a given fixed constant and 
the latter takes the difference between the current $J$ and the previous one, and terminates when the difference between the two falls below a certain prescribed tolerance.

Due to the nature of the iterative update presented in (\ref{equ_steepest_descent}), we expect each update on the control $\eta$ to reduce the objective functional $J$.
Hence we design an adaptive algorithm based upon this observation.
We may start with a arbitrary value of $\alpha$, namely $\alpha^\ell$.
If the computed objective function $J$ with $\alpha^\ell$ is smaller than the previous one, we increase the value of $\alpha^{\ell+1}$ and continue the computation.
However, if $J$ gets larger, this means the value $\alpha^\ell$ is not suitable, and the computation with $\eta^{\ell}$ need to be re-calculated using a new and smaller $\alpha$ (or the default minimum value $\alpha_{min}$).

We summarise this adaptive procedure in Algorithm \ref{alg_alpha}.
Note $\p_l$ and $\p_u$ are real numbers, unless otherwise stated, we set them to be $0.5$ and $1.1$, respectively.
In Section \ref{ada} we illustrate the effectiveness of this adaptive-$\alpha$ procedure.
\begin{algorithm}
\caption{Adaptive-$\alpha$}
\label{alg_alpha}
\begin{algorithmic}[0]
\State $\pmb{1. \text{While}}$ the difference between consecutive $J$s is still large or $J$ has not reached below a pre-defined tolerance $\pmb{\text{do}}$
\State $\pmb{2.}$ Solve the forward Allen-Cahn equation in $\Omega \times (0,T]$
\State $\pmb{3.}$ Compute the objective functional $J^\ell$
\State $\pmb{4.}$ $\pmb{\text{if}}$ $J^{\ell} > J^{\ell-1}$ $\pmb{\text{and}}$ $\ell > 0$ $\pmb{\text{then}}$
		\State \ \ \ \ \ \ \ \ \ \ \ $\alpha$ = $max(\alpha \times \p_l,\alpha_{min})$
		\State \ \ \ \ \ \ \ \ \ \ \ restart = $\pmb{\text{TRUE}}$
		\State \ \ \ \ $\pmb{\text{else if}}$ $J^{\ell} < J^{\ell-1}$ $\pmb{\text{and}}$ $\ell > 0$ $\pmb{\text{then}}$
		\State \ \ \ \ \ \ \ \ \ \ \ $\alpha$ = $\alpha \times \p_u$
		\State \ \ \ \ \ \ \ \ \ \ \ restart = $\pmb{\text{FALSE}}$
		\State \ \ \ \ $\pmb{\text{end if}}$
\State $\pmb{5.}$ $\pmb{\text{if}}$ restart == $\pmb{\text{FALSE}}$ $\pmb{\text{then}}$
		\State \ \ \ \ \ \ \ \ \ \ \ Solve the backward adjoint equation in $\Omega \times [T,0)$
		\State \ \ \ \ \ \ \ \ \ \ \ Backup the current $\eta$
		\State \ \ \ \ \ \ \ \ \ \ \ Compute the next $\eta$ using $\alpha$
		\State \ \ \ \ \ \ \ \ \ \ \ Continue to the next $\eta$ iteration
		\State \ \ \ \ $\pmb{\text{else}}$
		\State \ \ \ \ \ \ \ \ \ \ \ Compute a new $\eta$ using the latest backup with $\alpha$
		\State \ \ \ \ \ \ \ \ \ \ \ Restart the current $\eta$ iteration
		\State \ \ \ \ $\pmb{\text{end if}}$
\State $\pmb{6. \text{ End}}$
\end{algorithmic}
\end{algorithm}

\subsection{Space-time discretizations of the forward Allen-Cahn and adjoint equations}
\label{sec_32}

At the beginning of each $\eta$ iteration, we start by approximating the phase field equation (\ref{equ_Allen-Cahn}).
The spatial discretization scheme is a central finite difference method (FDM) with a standard seven-point stencil in 3-D on Cartesian grids with cell-centred vertices.
Although for illustrative purposes, the discrete system presented here is in 3-D, the use of a standard five-point stencil in 2-D is straightforward.
We assume $N$ is the number of grid points in each coordinate direction, $h$ is the uniform grid spacing (i.e. $h=\triangle x=\triangle y=\triangle z$), subscripts $i,j,k$ are used to indicate each grid point and each point has Cartesian coordinate $(x,y,z)$.
For the temporal discretization scheme we employ the fully-implicit second-order backward differentiation formula (BDF2) \cite{bdf2}.
With a given end time $T$, we assume a uniform time step size $\tau$.
For a time discrete sequence, $f$, we denote by $f^n:=f(t_n)$.
The standard BDF1 (also known as backward Euler method) is employed for the very first time step.

The result after applying the described discretisation to (\ref{equ_Allen-Cahn}) is the following algebraic system arising at each time step, find $\phi^{n+1,\ell+1}_{i,j,k}$ such that,
\begin{equation}
\label{equ_algebraic_Allen-Cahn}
\begin{split}
&\epsilon \frac{\phi^{n+1,\ell+1}_{i,j,k} - \frac{4}{3}\phi^{n,\ell+1}_{i,j,k} + \frac{1}{3}\phi^{n-1,\ell+1}_{i,j,k}}{\tau} =\\
& \frac{2\epsilon}{3} D\left(\phi_{i,j,k}^{n+1,\ell+1}\right) -\frac{2\left(-\phi^{n+1,\ell+1}_{i,j,k}+\left(\phi^{n+1,\ell+1}_{i,j,k}\right)^3\right)}{3\epsilon} +\frac{2\eta^{n+1,\ell+1}_{i,j,k}}{3}+\frac{2\lambda^{n+1}}{3},
\end{split}
\end{equation}
where $\ell+1$ denotes the current $\eta$ iteration, $n+1$, $n$ and $n-1$ indicate solutions from current, previous and the one before the previous time steps, respectively.
We denote the 3-D Laplacian operator $D$ as
\begin{equation}
\label{equ_laplacian}
D\left(\phi_{i,j,k}\right)= \frac{\phi_{i+1,j,k}+\phi_{i-1,j,k}+\phi_{i,j+1,k}+\phi_{i,j-1,k}+\phi_{i,j,k+1}+\phi_{i,j,k-1}-6\phi_{i,j,k}}{h^2}.
\end{equation}

Within each time step, while solving for the solution of the above system, we are also required to satisfy a given mass constraint.
This is done by iteratively determining the time-dependent, spatially-uniform volume constraint $\lambda$ for the imposed mass constraint \cite{Blowey}.
Therefore, the system in (\ref{equ_algebraic_Allen-Cahn}) has to be solved multiple times, until a stopping criterion for $\lambda$ is met.
We denote this $\lambda$ iteration using a superscript $\Lambda$, and its update follows the multi-step approach presented in \cite{Blowey}, which is given as
\begin{equation}
\label{equ_vc}
\lambda^{n+1,\Lambda+1} = \lambda^{n+1,\Lambda} + \frac{\left(\lambda^{n+1,\Lambda}-\lambda^{n+1,\Lambda-1}\right)\left[M_{\phi}^{n+1}-\int_{\Omega}\phi^{n+1,\Lambda}\right]}{\left(\int_{\Omega}\phi^{n+1,\Lambda}-\int_{\Omega}\phi^{n+1,\Lambda-1}\right)},\;\; \text{for $\Lambda>1$},
\end{equation}
where $M_\phi$ is defined in (\ref{equ_M}), $\Lambda+1$, $\Lambda$ and $\Lambda-1$ indicate values of $\lambda$ from current, previous and the one before the previous $\lambda$ iterations, respectively.
We follow \cite{Blowey} in using the initial guesses
\begin{equation}
\label{equ_lambda_guess}
\lambda^{\Lambda=0} = -\frac{2\epsilon}{\tau} + 1, \;\; \lambda^{\Lambda=1} = \frac{2\epsilon}{\tau} - 1.
\end{equation}
The stopping criterion used here is based upon the difference between consecutive values of $\lambda$.
Providing a tolerance $tol_\lambda$, we consider the algorithm to have converged when $|\lambda^{n+1,\Lambda+1}-\lambda^{n+1,\Lambda}| < tol_\lambda$.

From our experience, using the initial guesses in (\ref{equ_lambda_guess}) often led to more than three $\lambda$ iterations within each time step (with, say, a typical choice of $tol_\lambda = 0.01$).
For later computations, we can improve these initial guesses with known (already computed) values.
Within the first and second $\eta$ iterations, from time step $n=3$ onwards, we choose the computed $\lambda^{n-1, \ell}$ and $\lambda^{n, \ell}$, where $\ell = 1,2$, as our improved initial guesses.
When we are at the third $\eta$ iteration or beyond, we choose the two computed $\lambda$ (corresponding to the current time step) from the previous two $\eta$ iterations as initial guesses.
It must be observed that the computed solution $\phi$ at each time step has to be stored in order to compute the adjoint state $p$ later.

Having solved the algebraic system arising from the discretizations of the phase field representation of (\ref{equ_Allen-Cahn}) forward in time and stored the obtained solutions, we have only completed the first of two parts of the solution procedure.
The second part is to discretize and solve (\ref{equ_adjoint}).
We employ the described central FDM and BDF2 as our discretization schemes, and the resulting algebraic system for the adjoint state $p$ is the following,
\begin{equation}
\label{equ_algebraic_adjoint}
\begin{split}
&\frac{p^{n+1,\ell+1}_{i,j,k} - \frac{4}{3}p^{n+2,\ell+1}_{i,j,k} + \frac{1}{3}p^{n+3,\ell+1}_{i,j,k}}{\tau} =\\
&-\frac{2}{3}D\left(p_{i,j,k}^{n+1,\ell+1}\right) +\frac{2}{3}\left(\frac{-1+3\left(\phi^{n+1,\ell+1}_{i,j,k}\right)^2}{\epsilon^2}p^{n+1,\ell+1}_{i,j,k}\right).
\end{split}
\end{equation}
Note that the adjoint equation is posed backwards in time and its terminal condition is stated in (\ref{equ_adjoint}).
The BDF1 method is employed for the first time step.
For every subsequent time step the corresponding solution of $\phi$ that has been previously computed and stored enters as data in the adjoint equation.

\subsection{Techniques for improving algorithm efficiency}
\label{sec_33}

In this section, we explain several state-of-the-art algorithms that are used for obtaining the solution of (\ref{equ_algebraic_Allen-Cahn}) - (\ref{equ_algebraic_adjoint}) efficiently.
First of all, we describe an in-house two-grid solution strategy and this significantly improves the CPU time for solving the adjoint equation as well as massively reducing the memory requirement for storing all the solutions.
Secondly, we briefly mention a parallel, adaptive multigrid solution method that we use to solve the arising algebraic system at each time step.
Since this has been described in previous works, here we only briefly illustrate that the two-grid solution strategy can be a natural extension to the standard multigrid V-cycles.
We refer the interested reader to \cite{Bollada, Yang, Goodyer, Yang_paper} for further details.
Thirdly, the adaptive mesh refinement (AMR) can be included in the two-grid solution strategy in a straightforward manner.
This significantly improves the CPU time for solving the Allen-Cahn equation.
The final technique is the adaptive-$\alpha$ algorithm which has already been mentioned in Section \ref{adaptive-a}, here we omit details of this algorithm.

We propose a two-grid solution strategy that exploits a key difference between the forward Allen-Cahn and the backward adjoint equations.
As is well known for the Allen-Cahn equation the parameter $\epsilon$ determines the thickness of the diffuse interfacial region, $\Gamma_\epsilon$, that approximates the hypersurface $\Gamma$.
In order for the Allen-Cahn equation to reliably approximate mean curvature flow the interfacial region has to be well resolved.
Typically $8$ grid points are required across the width of the diffusive interface, see \cite{DDE}.
On the other hand, in our numerical simulations we observed that the solution of the backward adjoint equation varies less, see Figures \ref{fig_four} (c) and \ref{fig_four} (d) in which the optimised solutions of the adjoint equation at time $t=0.0625$ and $t=0.125$ are displayed, and so a milder restriction on the grid size $h$ in the interfacial region is expected when solving this equation.
Our numerical tests (see Section \ref{sec_results}) suggest that such a strategy can dramatically reduce the memory requirement and increase computational efficiency without significantly compromising accuracy.

\begin{figure}[!htbp]
\centering
\fbox{
\includegraphics[width=16cm]{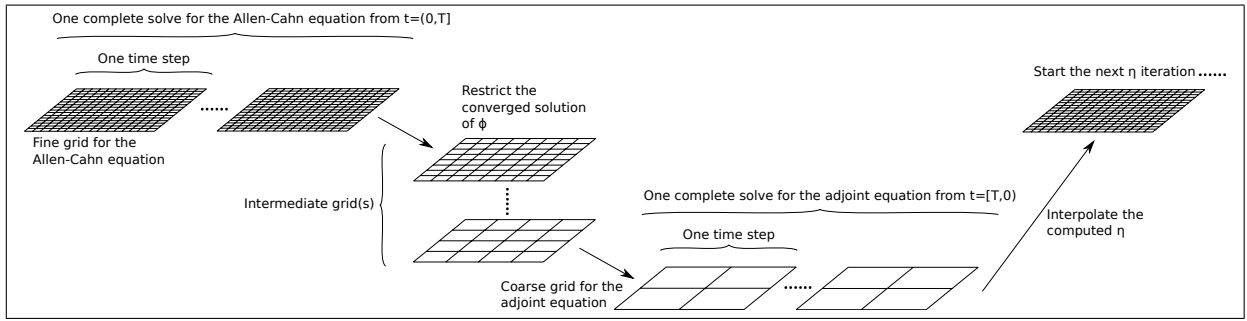}
}
\caption{Sketch illustrating our in-house two-grid solution strategy, where the adjoint equation is solved on a much coarser grid. The storage for all-time-step solutions is done on such a grid so as to reduce the memory requirement.}
\label{fig_2_two_grid}
\end{figure}

Within this implemented robust in-house two-grid solution strategy, we solve the Allen-Cahn equation on a grid hierarchy where its finest grid has sufficient grid points for the chosen $\epsilon$; the backward adjoint equation is then solved using only part of the grid hierarchy to improve the efficiency.
As we stated previously, during the computation, at least solutions from two variables of $\phi$, $p$ and $\eta$ from all time steps in the current $\eta$ iteration are stored.
If a very fine mesh resolution is required, this undoubtedly imposes a severe requirement for the memory, even in a parallel setting.
To relax this constraint, we store all-time-step solutions only on the coarser grid where we solve the adjoint problem.
When the stored solutions are required on finer grids (where $\phi$ is solved), an interpolation is used to transfer the stored solutions to that grid level.
We illustrate this two-grid solution strategy in Figure \ref{fig_2_two_grid}.

\begin{remark}
\label{rem1}
When solving such optimal control problems we note that apart from the computational complexity, in terms of CPU time, it must be noted that there is also a very large memory requirement.
More specifically, $\eta$ exists on all internal grid points for all time steps.
The $\eta$ update \eqref{equ_steepest_descent} requires the solution of the adjoint equation $p$ on every grid point at every time step.
For the adjoint equation \eqref{equ_adjoint}, the solution of the state equation $\phi$ on every space-time grid point is also required.
Thus, when we solve the Allen-Cahn equation, the solutions of $\phi$ and $\eta$ are stored on every grid point for all time steps.
The requirements on memory storage can become significant as the number of grid point increases.

To give the reader a brief idea, say one uses a 2-D grid consisting of $512^2$ grid points with $50$ forward and backward time steps for the Allen-Cahn and the adjoint equation respectively.
The solutions of $2$ variables ($\phi$ and $\eta$) are stored using a double-precision format, i.e., each float-point number occupies $8$ bytes in the memory.
This setting requires around $210$ megabytes of memory space.
However, if this simulation is done in 3-D on a uniform grid with the resolution of $512^3$, the memory requirement becomes approximately $107$ gigabytes.
Such a simulation may not be feasible if no parallelism nor any other advanced techniques are employed.
To this end, we couple our two-grid solution strategy with a parallelisation of a domain decomposition approach, as well as dynamic AMR to enable 3-D simulations.
\end{remark}

\begin{figure}[!htbp]
\centering
\fbox{
\includegraphics[width=15cm]{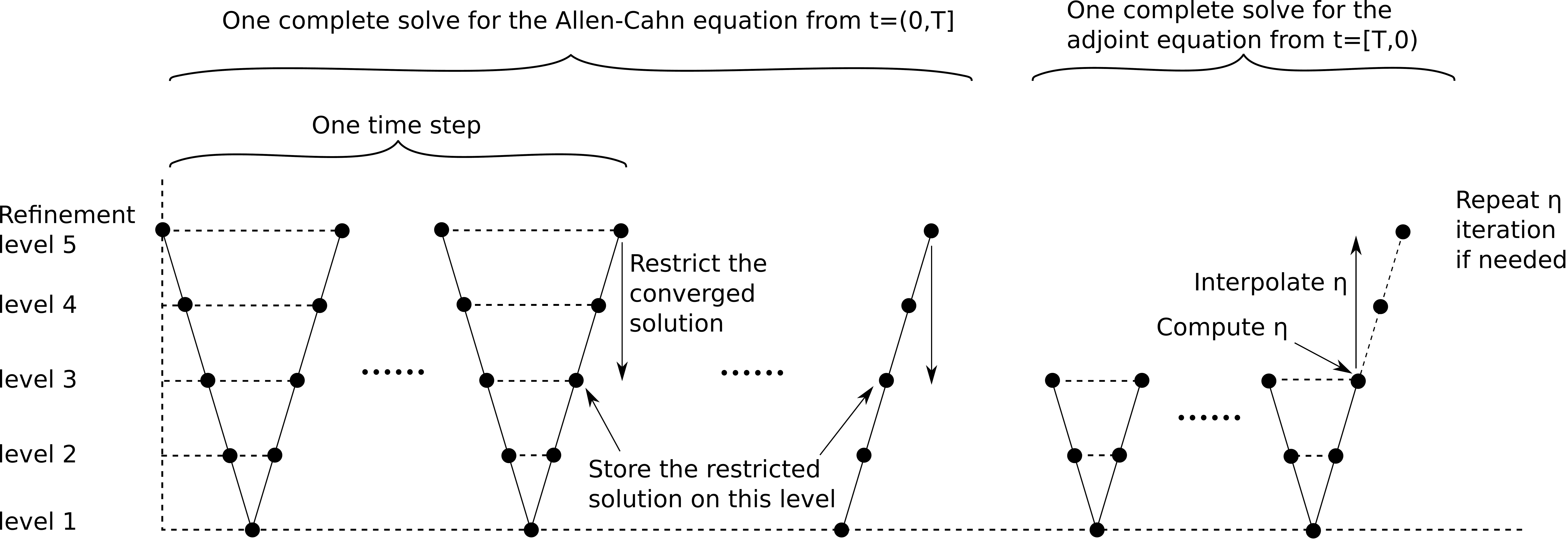}
}
\caption{Sketch demonstrating our in-house multi-depth V-cycle multigrid strategy where the adjoint equation is solved on a much coarser grid. The storage for all-time-step solutions is done on such a grid so as to reduce the memory requirement.}
\label{fig_2_multi-depth}
\end{figure}

The software framework used here is called Campfire v2.0.
Comparing with its previous versions used in \cite{Bollada, Yang, Yang_paper}, the latest software received some significant changes to its structure in order to deal with a forward and a backward solve, as well as additional parallel memory allocations.
This framework contains a geometric nonlinear multigrid solution method with a full approximation scheme (FAS), as well as a multi-level adaptive technique (MLAT) variant.
However, when solving the linear adjoint equation, this FAS multigrid reduces to the standard linear multigrid method \cite{bigMG:2001}.
The multigrid methods are widely known to be one of the fastest numerical methods with a linear complexity \cite{bigMG:2001,Brandt,Briggs:2000}, and we demonstrate this later in the paper.
Its parallelisation comes from a domain decomposition technique, and message passing interface is used for parallel communication.
Here using Figure \ref{fig_2_multi-depth}, we illustrate that the presented in-house two-grid solution strategy can be comfortably extended into the multigrid V-cycles.

Considering we are using phase field approximation, dynamic AMR is commonly employed so the interfacial region is well captured and resolved and the computations are saved in other regions where the phase field variable tends to be a constant value.
On the other hand, if the adjoint equation is solved with the dynamic AMR, it would involve further complications where we need to store the evolved mesh structure.
This can be practically unfeasible with the memory requirement that is mentioned in the Remark \ref{rem1} earlier.
Therefore, we employ dynamic AMR for the solutions of the forward Allen-Cahn equationand solve the backwards adjoint PDE on a fixed (and much coarser) grid.
For instance, the finest grid shown in Figure \ref{fig_2_two_grid} may be a dynamically adaptive grid, as well as the intermediate grids.
However, the coarse grid that we solve the adjoint equation stays fixed.
We refer the interested reader to \cite{Yang, Yang_paper} for the description of the dynamic AMR and its associated dynamic load-balancing in parallel.

\section{Numerical experiments}
\label{sec_results}

All the results shown in this section were generated using the local HPC cluster provided and managed by the University of Sussex.
This HPC cluster consists of $3000$ computational units.
The models of the computational units are AMD64, x86\_64 or 64 bit architecture, made up of a mixture of Intel and AMD nodes varying from 8 cores up to 64 cores per node.
Each unit is associated with $2$GB memory space.
Most of the simulations in this paper were executed using $4-32$ cores.
The parallel scalability of our multigrid solver, Campfire, has been discussed in earlier publications such as \cite{Bollada, Yang, Yang_paper}, where in \cite{Bollada}, they successfully scaled up to one thousand computational cores on the national supercomputer VECToR in 2013.
We refer the reader to \cite{Yang, Yang_paper} for a detailed explanation and related results for the parallel scalability of our software.

\subsection{A 2-D benchmark example}
\label{ss4p1}
\begin{figure}[!htbp]
\centering
\includegraphics[width=14.5cm]{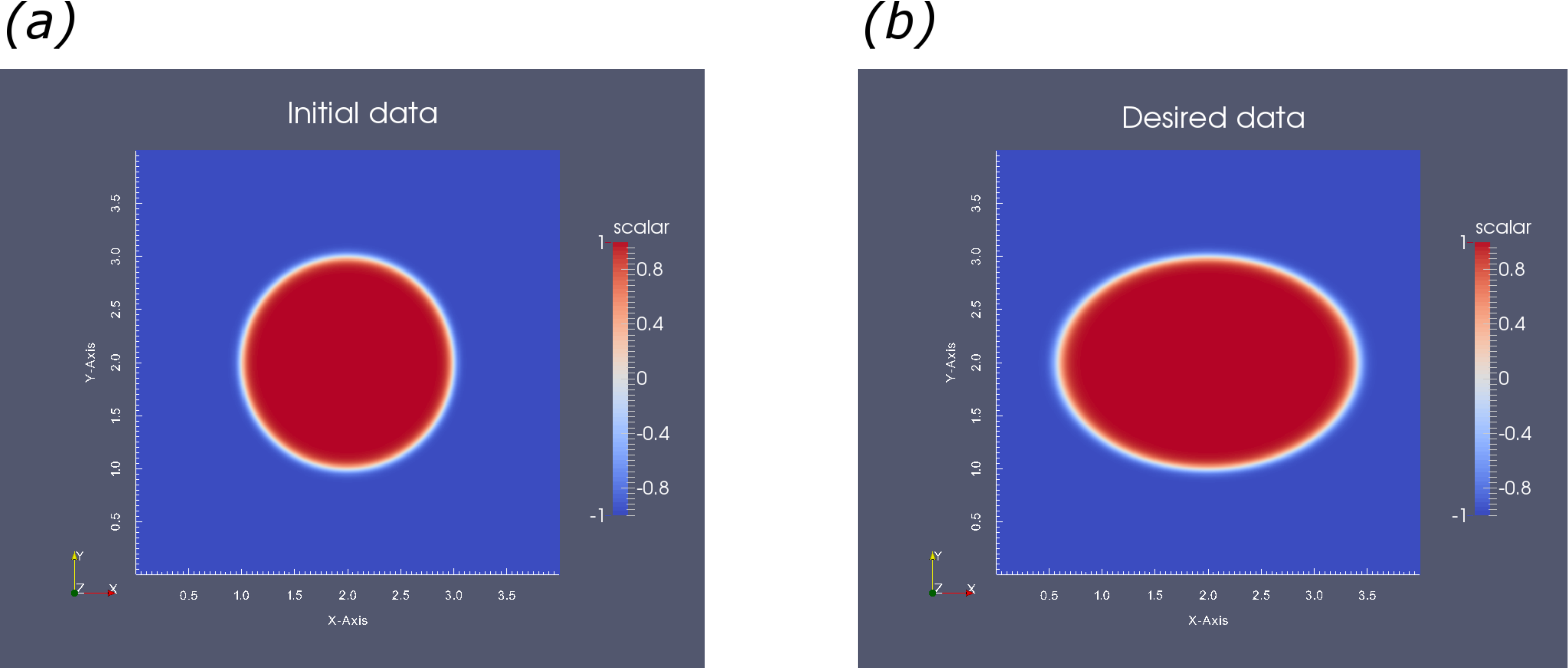}
\caption{(a) shows the initial data (i.e. (\ref{equ_2_initial})) and (b) shows the desired data (i.e. (\ref{equ_2_desired})). The colour version of this figure is online.}
\label{fig_ini_de}
\end{figure}

We start with a benchmark 2-D example.
The initial data is a circle centred at $(2,2)$ with radius $1$.
We use a hyperbolic tangent function to obtain a continuous interfacial region with a width of $\mathcal{O}(\epsilon)$ 
\begin{equation}
\label{equ_2_initial}
\phi^{t=0} = \tanh\left( \frac{-\left[ \left(x-2\right)^2+\left(y-2\right)^2-1\right]}{\epsilon} \right).
\end{equation}
The desired data is an ellipse:
\begin{equation}
\label{equ_2_desired}
\phi_{obs} = \tanh\left( \frac{-\left[ \frac{\left(x-2\right)^2}{2}+\left(y-2\right)^2-1\right]}{\epsilon} \right).
\end{equation}
Both the initial and desired shapes are illustrated in Figure \ref{fig_ini_de}.

For illustrative purposes, we take the computational domain $\Omega=(0,4)^2$.
The choices of parameters are given in Table \ref{table_parameters}.
At each time step, the infinity norm of the residual is assessed, and the computation is said to be converged when this norm is smaller than $1.0\times 10^{-11}$.
\begin{table}[!htb]
\begin{center}
\begin{tabular}{|c|c|c|}
\hline
Name& Description &Value \\
\hline
$\alpha$			 & Step size for the control update& $0.1$\\
$\theta$ 			 & Regularisation parameter        &$0.01$\\
$\epsilon$ 		 & Width of the diffuse interface  & $0.1$  \\
$T$			 			 & Total time duration    				 &$0.125$ \\
$\tau$ 				 & Time step size 								 &(varies)\\
\hline
\end{tabular}
\end{center}
\caption{The parameters of the optimal control problem for the 2-D examples.}
\label{table_parameters}
\end{table}

We present the computed solutions from a uniform grid with a resolution of $1024^2$.
We select a time step $\tau = 7.8125\times 10^{-4}$, which yields $160$ time steps.
For the propose of demonstration, we run $50$ $\eta$ iterations. 
In a multigrid setting, we use a $16^2$ as the coarsest grid, and grids like $32^2$, $64^2$ etc. are the intermediate grids in the V-cycle hierarchy.
The multigrid hierarchy is used for all simulations ($16^3$ used as the coarsest grid for 3-D simulations) and is fairly standard, for clarity we forgo mentioning the multigrid setting later on.

\begin{figure}[!htbp]
\centering
\includegraphics[width=12cm]{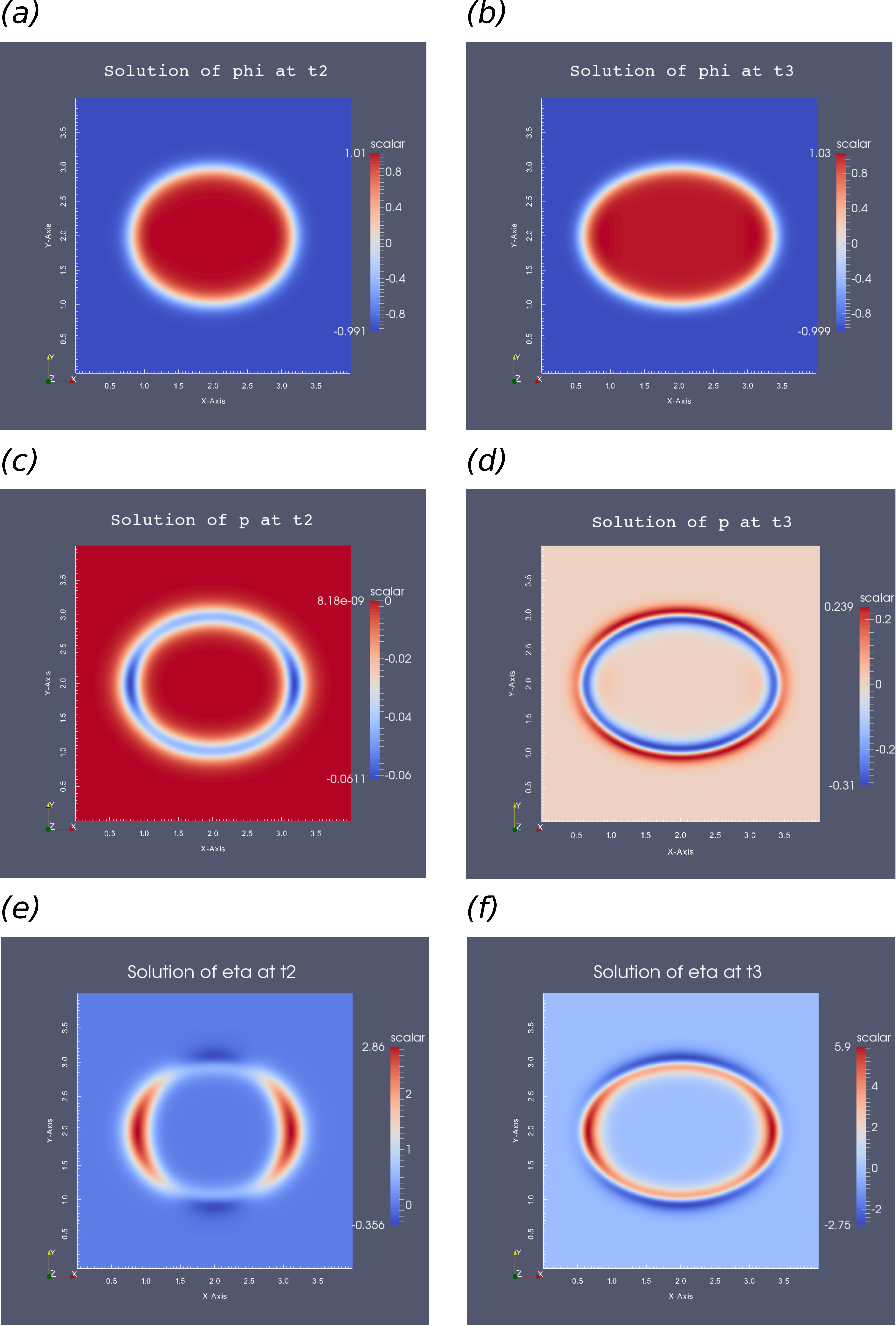}
\caption{(a) illustrates the computed solution of $\phi$ halfway through (i.e. $t=0.0625$) and (b) shows the computed final shape (i.e. $t=T=0.125$. The corresponding solutions of adjoint $p$ are included in (c) and (d) and the solutions of $\eta$ are shown in (e) and (f). The colour version of this figure is online.}
\label{fig_four}
\end{figure}

The computed solution of $\phi$ at $t=0.0625$ (halfway through the time series) is presented in Figure \ref{fig_four} (a).
The computed final shape after the last time step is illustrated in Figure \ref{fig_four} (b).
The corresponding solutions of the adjoint $p$ at time $t=0.0625$ and $t=T=0.125$ are included in Figures \ref{fig_four} (c) and (d).
The corresponding solutions of $\eta$ are shown in Figures \ref{fig_four} (e) and (f).
We can see that, as expected for a circle evolving into an ellipse, the most forcing in (c) is placed at the left and the right.
We also observe that in (d), the forcing is highly positive on the inner side of the phase field interface, and highly negative on the other, in order to keep the shape from further unwanted expansion or shrinking.
The reductions in the objective functional $J$ from the $50$ $\eta$ iterations are shown using a semi-log plot in Figure \ref{fig_adaptive_alpha}.

\subsection{Convergence tests for the benchmark 2-D example}
\label{conv}
\begin{table}[!htb]
\begin{center}
\begin{tabular}{|c|c|c|c|}
\hline
\multicolumn{4}{|c|}{$\text{$L_2(\Omega)$ error for $\phi$}$} \\
\hline
  & $m = t_1$ & $m = t_2$ & $m=t_3$\\
\hline
$d^m_{64^2}$  & $3.4264 \times 10^{-2}$ & $4.9226\times 10^{-2}$ & $6.8561 \times 10^{-2}$\\
$d^m_{128^2}$ & $1.9721 \times 10^{-2}$ & $4.8058\times 10^{-2}$ & $6.0397 \times 10^{-2}$\\
$d^m_{256^2}$ & $8.1793 \times 10^{-3}$ & $2.2300\times 10^{-2}$ & $3.4557 \times 10^{-2}$\\
$d^m_{512^2}$ & $2.7850 \times 10^{-3}$ & $8.0407\times 10^{-3}$ & $1.3192 \times 10^{-2}$\\
\hline
\multicolumn{4}{|c|}{$\text{$L_2(\Omega)$ error for adjoint $p$}$} \\
\hline
  & $m = t_1$ & $m = t_2$ & $m=t_3$\\
\hline
$d^m_{64^2}$  & $1.6773 \times 10^{-2}$ & $1.8048\times 10^{-2}$ & $4.9344 \times 10^{-2}$\\
$d^m_{128^2}$ & $9.9721 \times 10^{-3}$ & $1.0158\times 10^{-2}$ & $3.1554 \times 10^{-2}$\\
$d^m_{256^2}$ & $7.9290 \times 10^{-3}$ & $8.5311\times 10^{-3}$ & $2.2551 \times 10^{-2}$\\
$d^m_{512^2}$ & $6.5082 \times 10^{-3}$ & $7.5551\times 10^{-3}$ & $1.4901 \times 10^{-2}$\\
\hline
\multicolumn{4}{|c|}{$\text{$L_2(\Omega)$ error for $\eta$}$} \\
\hline
  & $m = t_1$ & $m = t_2$ & $m=t_3$\\
\hline
$d^m_{64^2}$  & $1.6976 \times 10^{-1}$ & $2.0752\times 10^{-1}$ & $7.5240 \times 10^{-1}$\\
$d^m_{128^2}$ & $1.1923 \times 10^{-1}$ & $1.5793\times 10^{-1}$ & $6.2554 \times 10^{-1}$\\
$d^m_{256^2}$ & $8.5093 \times 10^{-2}$ & $1.0601\times 10^{-1}$ & $5.3023 \times 10^{-1}$\\
$d^m_{512^2}$ & $3.2344 \times 10^{-2}$ & $3.9359\times 10^{-2}$ & $2.6302 \times 10^{-1}$\\
\hline
\end{tabular}
\end{center}
\caption{The convergence tests for the solutions of $\phi$, adjoint $p$ and $\eta$.}
\label{table_eta}
\end{table}

In this subsection, we report on numerical evidence that the proposed optimal control model converges as we refine both spatially and temporally.

We use the simulation that was described in the previous subsection as a benchmark.
To recap, it was solved on a grid with the resolution of $1024^2$ and $160$ number of time steps.

In order to conduct the convergence tests, the optimal control model here is solved independently on the following grids: $64^2$, $128^2$, $256^2$ and $512^2$.
We use a time step size $\tau = 0.0125$ for the $64^2$ simulation, which has $10$ time steps.
Then the choices of $\tau$ is halved each time we use a finer grid.
This results in $128^2$ to have $20$, $256^2$ to have $40$ and $512^2$ to have $80$ time steps, respectively.
The parameters for these simulations are the same as shown in Table \ref{table_parameters}.
However, instead of running the $\eta$ iteration to a constant, we change the stopping criterion so that the $\eta$ iteration stops when the objective functional $J$ is below $0.065$.
Note this new stopping criterion requires 24 $\eta$ iterations to be satisfied for the $1024^2$ simulation and roughly the same for other simulations.

We assess the solutions generated at three specific times.
They are $t_1 = 0.0125$, $t_2 = T/2 = 0.0625$ and $t_3 = T = 0.125$.
Note $t_1$ is at the end of first time step of the $64^2$ simulation, and subsequently the end of second, fourth and eighth time steps of the $128^2$, $256^2$ and $512^2$ simulations respectively.

To compare the solutions spatially at the corresponding $t_1$, $t_2$ and $t_3$, further procedures are required.
This is because the solutions are generated with different resolutions of grids.
As mentioned earlier, we use the solution from $1024^2$ grid as the benchmark, therefore all the solutions from other simulations other than $1024^2$ are interpolated to the uniform $1024^2$ grid.
Note this process uses a standard bilinear interpolation \cite{bigMG:2001}, which is also the one used in our multigrid solver.

The interpolated solutions are compared with the benchmark solutions at $t_1$, $t_2$ and $t_3$.
Solutions from all grid points are assessed for the difference from the benchmark solution in the $L_2(\Omega)$ error defined as follows
\begin{equation}
\label{equ_2_norm}
d^m_l := \frac{\sum_{i=1}^{N}\sum_{j=1}^{N}(\phi^{m,l}_{i,j} - \phi_{i,j}^{m,1024^2})^2}{N\times N},~~~m=t_1,t_2,t_3,~l=64^2,128^2,256^2,512^2
\end{equation}
where $\phi^{m,l}_{i,j}$ is the computed and interpolated value of $\phi$ on grid $l$ and $N=1024$ is the number of internal grid points (after interpolation) on one axis.
This is repeated for the solutions of the adjoint $p$ and $\eta$.
We summarise the convergence tests in Table \ref{table_eta}.
It can be seen from this table that as we refine both spatially and temporally, the solutions of the proposed optimal control model appear to converge.

\subsection{The multigrid performance on the benchmark 2-D example}

Within our software framework, the algebraic system arising from each time step is solved by a multigrid solver.
Here in this subsection, we assess, numerically, the performance of our multigrid performance.

\begin{figure}[!htbp]
\centering
\includegraphics[width=15cm]{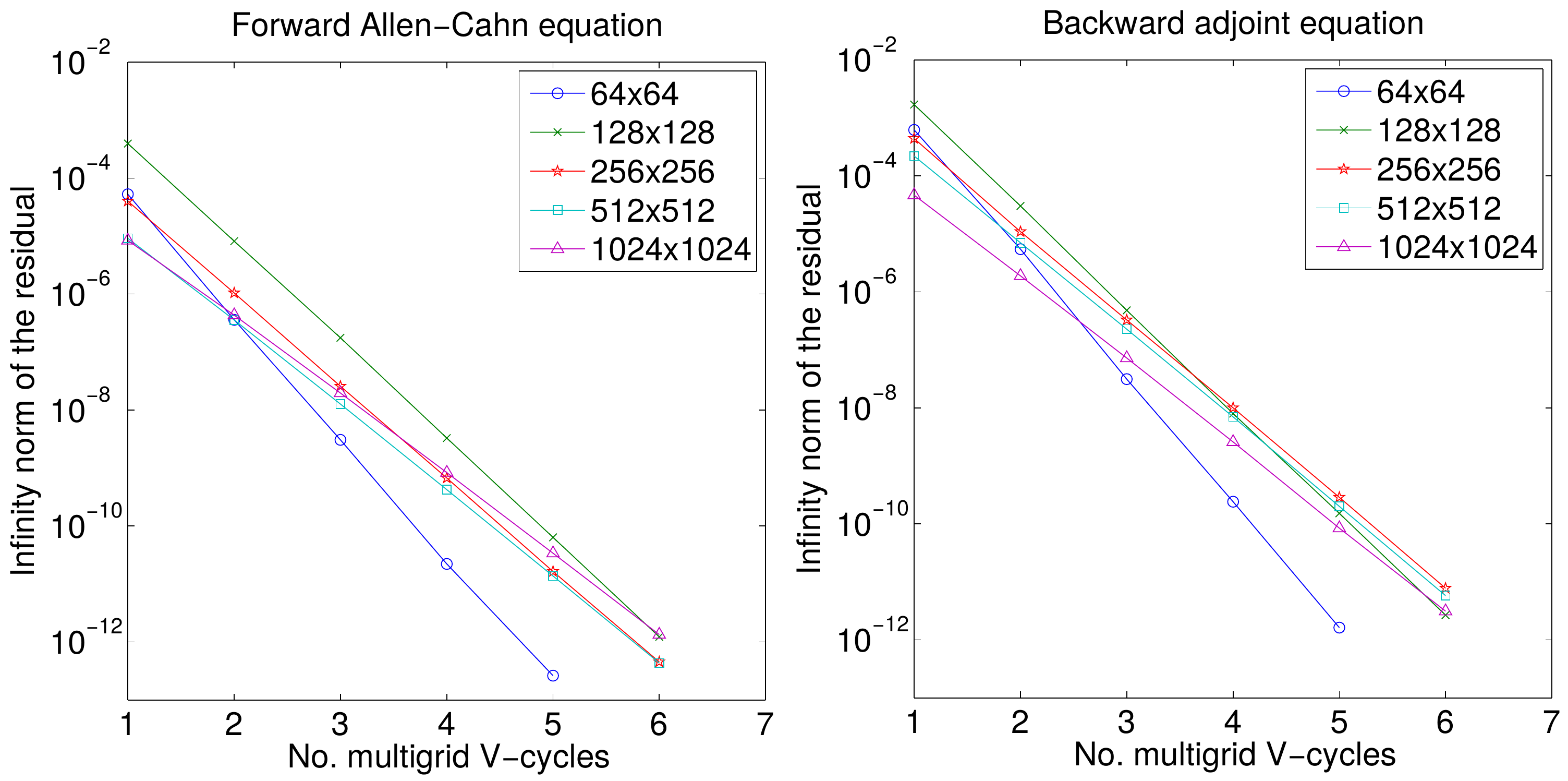}
\caption{The multigrid convergence rates for the forward Allen-Cahn and backward adjoint equations. The colour version of this figure is online.}
\label{fig_2D_V-cycle}
\end{figure}

First of all, we present multigrid convergence rates by plotting the infinity norm of the residual at the end of each V-cycle from a typical time step.
Furthermore, since two different equations are solved separately, we separate them and illustrate these results in Figure \ref{fig_2D_V-cycle}.
All eight lines from two plots in Figure \ref{fig_2D_V-cycle} are nearly parallel to each other which suggests the reductions in the infinity norms of the residuals are independent of grid sizes.

\begin{figure}[!htbp]
\centering
\includegraphics[width=8cm]{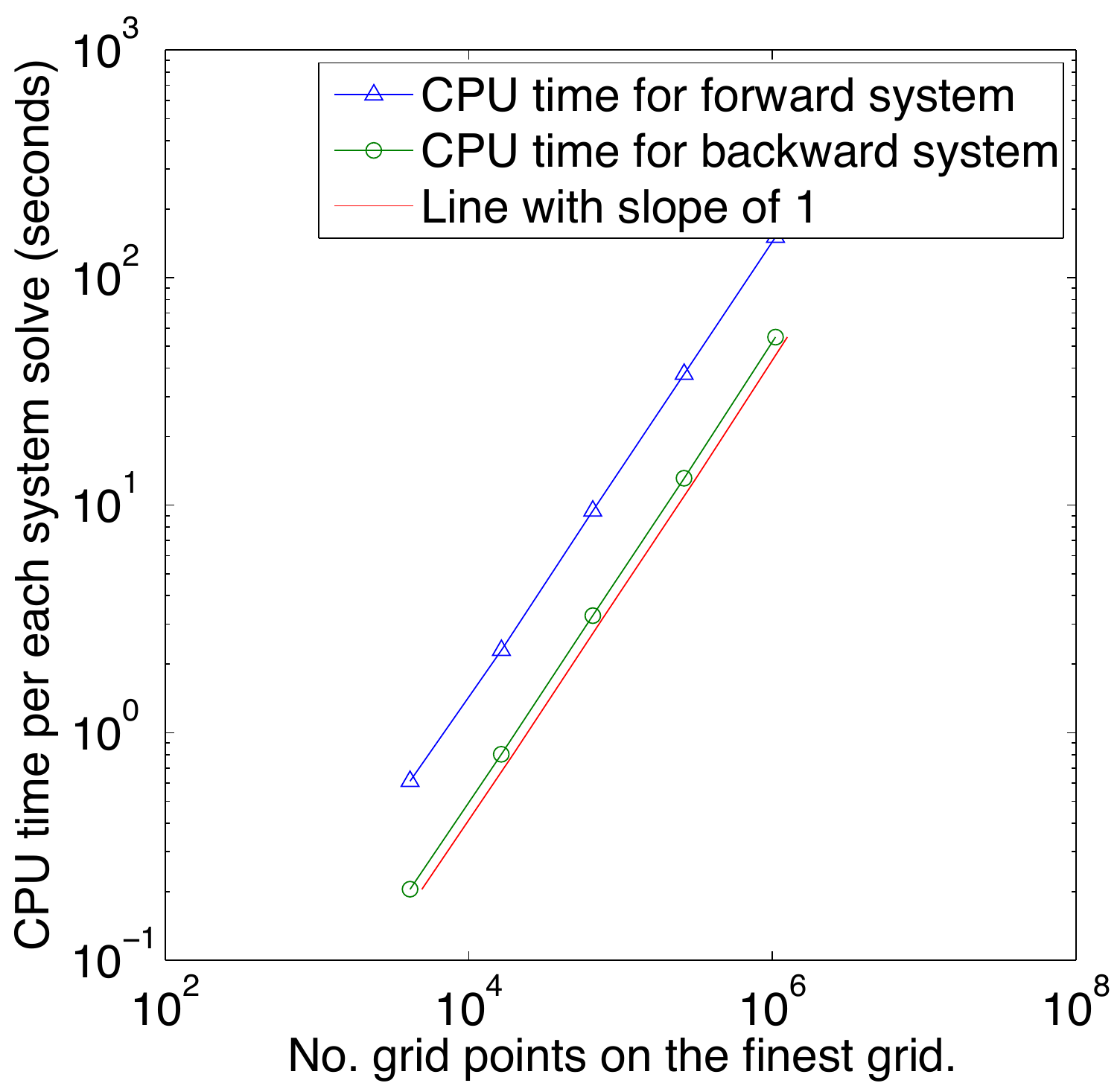}
\caption{A log-log plot to illustrate the linear complexity of our multigrid solver. For comparisons, a line of slop $1$ is included. The colour version of this figure is online.}
\label{fig_time}
\end{figure}

We demonstrate the linear complexity of our multigrid solver in Figure \ref{fig_time}.
Five simulations ($64^2$, $128^2$, $256^2$, $512^2$ and $1024^2$) are timed with a single computational core and the averages of the CPU costs for $10$ $\eta$ iterations are plotted as shown in Figure \ref{fig_time}.
For clarity, we plot a red line of slope $1$.
The reason why the forward solver costs more is because, the mass constraint typically requires $2$ - $3$ choices of $\lambda$ (more in the first two $\eta$ iterations, since the guesses for $\lambda$ are poor), meaning the algebraic system from the forward equation has to be solved $2$ - $3$ times within each time step.
On the other hand, the algebraic system for the backward equation only requires to be solved once.

\subsection{Adaptive-$\alpha$ algorithm with the benchmark 2-D example}
\label{ada}

An additional improvement to the efficiency may come from using the described adaptive-$\alpha$ algorithm (see Algorithm \ref{alg_alpha}).
In this subsection, we illustrate the effectiveness of this approach.

\begin{figure}[!htbp]
\centering
\includegraphics[width=15cm]{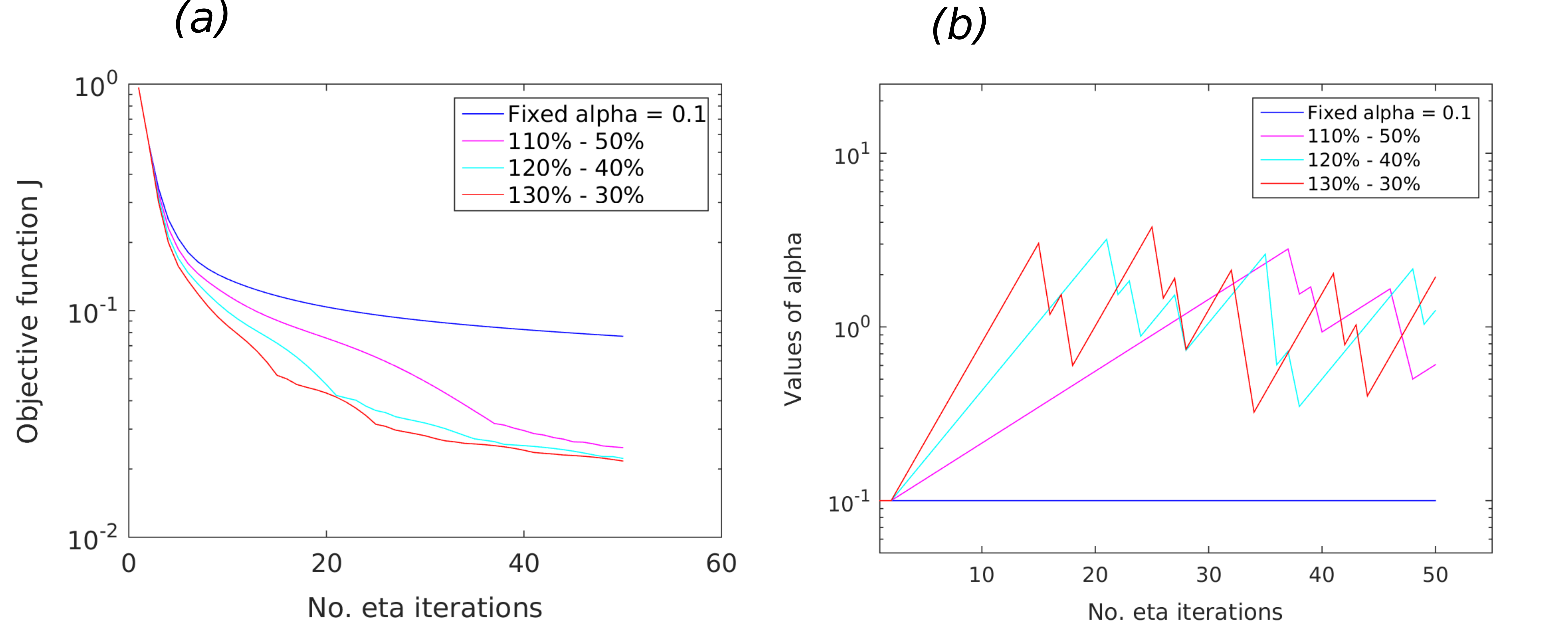}
\caption{A semi-log plot shows the reductions of the objective functional $J$ from using constant and adaptive $\alpha$s in (a). A semi-log plot shows the changes in the values of $\alpha$ in (b). The colour version of this figure is online.}
\label{fig_adaptive_alpha}
\end{figure}

The two most influential parameters in the algorithm are $\p_l$ and $\p_u$ which control the incremental and decremental portions of the step size respectively.
For the purpose of demonstration, we use the described $1024^2$ simulation with $160$ time steps.
As mentioned previously, the reductions of the objective functional $J$ within a fixed $50$ $\eta$ iterations is shown in Figure \ref{fig_adaptive_alpha} (a).
The initial choice of $\alpha^{\ell=0} = 0.1$ is the same used in the fixed $\alpha$ simulations previously.

Here we choose three different incremental parameters: $\p_u^1 = 110\%$, $\p_u^2 = 120\%$ and $\p_u^3 = 130\%$.
The corresponding decremental parameters are $\p_l^1 = 50\%$, $\p_l^2 = 40\%$ and $\p_l^3 = 30\%$, respectively.
Thus the more increases to the step size, the harder the penalisations become.
For comparison, simulations with these three pairs of parameters are done with $50$ fixed $\eta$ iterations.
We plot the reductions of the objective functional $J$ from using the adaptive-$\alpha$ algorithm in Figure \ref{fig_adaptive_alpha} (a).
Note that in this figure, we only show the $J$ from the successful iterations.
For completeness, we illustrate the evolutions of $\alpha$ from all four simulations in Figure \ref{fig_adaptive_alpha} (b), where the decreases in the values of $\alpha$ indicate the failed attempts.
A trade-off can be observed from these two figures: a larger incremental parameter leads to a faster convergence, however, this may result in more failed attempts and thus in turn results in more computational time.
In this case the gains in efficiency of the adaptive-$\alpha$ approach against a fixed value of $\alpha$ are evident. The use of an adaptive $\alpha$ is motivated by the fact that in general our initial guess for the solution to the optimal control problem may be poor and hence large step sizes may be admissible in the steepest descent update as we are far from local minima. As we approach the local minima smaller step sizes are necessary to prevent overshoot and hence some adaptivity in the parameter $\alpha$ is expected to be desirable.
More involved algorithms for the selection of an optimal parameter $\alpha$, e.g., via a line search \cite{Troltzsch}, may also be worthwhile topics of future investigation.

\subsection{Two-grid solution strategy with dynamic AMR on the benchmark 2-D example}

In Remark \ref{rem1} we note that, due to the very large number of degrees of freedom that are typically required in order to accurately resolve  phase field representation of interfaces with large curvatures, simulations on uniform grids in three space dimensions are unlikely to be feasible. To this end we introduce the two-grid AMR solution strategy proposed in Section \ref{sec_33}.
In this subsection, we investigate the effectiveness of this two-grid solution strategy, in two space dimensions, and illustrate its robustness with the use of dynamic AMR.

We consider a two-grid simulation where we solve the forward Allen-Cahn equation on a $1024^2$ uniform grid 
while the adjoint equation and the storage of all the solutions, $\eta$, $\phi$, $p$, takes place on a $64^2$ uniform grid.
The solutions of this two-grid simulation are compared with solutions using a standard (one-grid) $64^2$ uniform grid simulation.
We note that in both simulations all solutions are stored on a grid with resolution of $64^2$.

For simplicity, we take $160$ time steps for both simulations so that temporal errors have less influence.
Like the convergence tests shown in Subsection \ref{conv}, we solve the system until $J$ gets below $0.065$.
In order to compute the error, the solutions from both simulations are interpolated and compared against solutions from the $1024^2$ simulation.

We illustrate the errors in Table \ref{table_two_grid_1}, where 
$$
d^m_{1024^2-64^2} := \frac{\sum_{i=1}^{N}\sum_{j=1}^{N}(\phi^{m,1024^2-64^2}_{i,j} - \phi_{i,j}^{m,1024^2})^2}{N\times N},~~~m=t_1,t_2,t_3, 
$$
with $\phi^{m,1024^2-64^2}_{i,j}$ denoting the $\phi$ solution from the two-grid simulation.

From this table we see that solving $\phi$ on a finer grid, while solving for $p$ on a coarse grid and storing all solutions on the coarse grid, 
not only results in a reduction of the error in $\phi$ but it also results in a reduction of the errors of $p$ and $\eta$.
\begin{table}[!htb]
\begin{center}
\begin{tabular}{|c|c|c|c|}
\hline
\multicolumn{4}{|c|}{$\text{$L_2(\Omega)$ error for $\phi$}$} \\
\hline
  & $m = t_1$ & $m = t_2$ & $m=t_3$\\
\hline
$d^m_{64^2}$        & $3.1224 \times 10^{-2}$ & $4.7216\times 10^{-2}$ & $6.3531 \times 10^{-2}$\\
$d^m_{1024^2-64^2}$ & $8.1840 \times 10^{-3}$ & $1.6194\times 10^{-2}$ & $8.6616 \times 10^{-3}$\\
\hline
\multicolumn{4}{|c|}{$\text{$L_2(\Omega)$ error for adjoint $p$}$} \\
\hline
  & $m = t_1$ & $m = t_2$ & $m=t_3$\\
\hline
$d^m_{64^2}$        & $1.3722 \times 10^{-2}$ & $1.2018\times 10^{-2}$ & $4.0874 \times 10^{-2}$\\
$d^m_{1024^2-64^2}$ & $6.1444 \times 10^{-3}$ & $4.9266\times 10^{-3}$ & $9.9998 \times 10^{-3}$\\
\hline
\multicolumn{4}{|c|}{$\text{$L_2(\Omega)$ error for $\eta$}$} \\
\hline
  & $m = t_1$ & $m = t_2$ & $m=t_3$\\
\hline
$d^m_{64^2}$        & $1.2176 \times 10^{-1}$ & $1.8872\times 10^{-1}$ & $6.9943 \times 10^{-1}$\\
$d^m_{1024^2-64^2}$ & $9.9110 \times 10^{-2}$ & $1.3768\times 10^{-1}$ & $4.2063 \times 10^{-1}$\\
\hline
\end{tabular}
\end{center}
\caption{Comparisons of errors between a two-grid simulation ($1024^2-64^2$) and a standard $64^2$ simulation.}
\label{table_two_grid_1}
\end{table}
However since the number of degrees of freedom in the two-grid ($1024^2-64^2$) simulation is considerably larger than the number of degrees of freedom in the 
standard $64^2$ simulation. This behaviour is somewhat expected. In the next simulation we conduct a comparison between a two-grid simulation that solves the Allen-Cahn equation on an adaptive $256^2$ grid with dynamic AMR and the adjoint equation on a uniform $64^2$ grid, with a simulation on a standard $128^2$ uniform grid. In this comparison the number of degrees of freedom in the two simulations is comparable, $17200$ (maximum number of degrees of freedom occurred)  for the two-grid simulation versus $128^2=16384$.

\begin{figure}[!htbp]
\centering
\includegraphics[width=13cm]{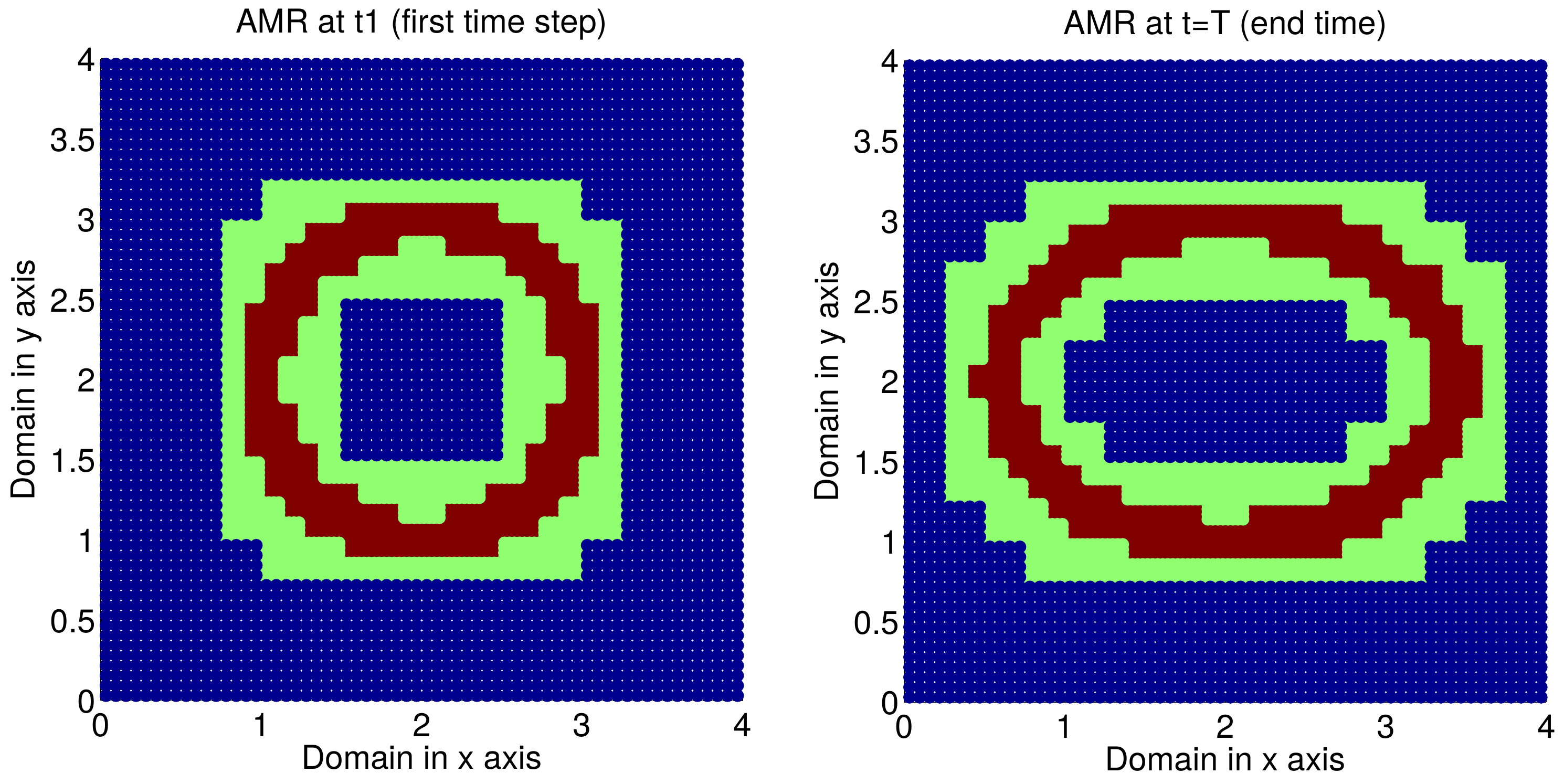}
\caption{Two colour plots show the dynamic AMR in our solver. The blue region shows the $64^4$ grid; light green region indicates the $128^2$ grid; and finally red region illustrates the finest $256^2$ grid. The colour version of this figure is online.}
\label{fig_multi-depth_AMR}
\end{figure}

Both simulations have $160$ time steps. 
The errors are shown in Table \ref{table_two_grid_2}.
From this table, we can see that for the two-grid simulation only the errors in $\phi$ are better.
This is expected as the adjoint is solved on a coarser grid (i.e. $64^2$) and $\eta$, $\phi$ and $p$ are stored on this coarse grid. 
On the other hand, it is important to note that we can store all the solutions on a coarser grid as well as solving the adjoint equation there without compromising too much on accuracy; this is crucial for 3-D simulations. 
In Figure \ref{fig_multi-depth_AMR} we show two snapshots of our dynamic AMR at $t=t_1$ and $t=T$.
\begin{table}[!htb]
\begin{center}
\begin{tabular}{|c|c|c|c|}
\hline
\multicolumn{4}{|c|}{$\text{$L_2(\Omega)$ error for $\phi$}$} \\
\hline
  & $m = t_1$ & $m = t_2$ & $m=t_3$\\
\hline
$d^m_{128^2}$      & $1.2327 \times 10^{-2}$ & $2.6664\times 10^{-2}$ & $3.8450 \times 10^{-2}$\\
$d^m_{256^2-64^2}$ & $8.6270 \times 10^{-3}$ & $1.6895\times 10^{-2}$ & $3.2925 \times 10^{-2}$\\
\hline
\multicolumn{4}{|c|}{$\text{$L_2(\Omega)$ error for adjoint $p$}$} \\
\hline
  & $m = t_1$ & $m = t_2$ & $m=t_3$\\
\hline
$d^m_{128^2}$      & $9.2021 \times 10^{-3}$ & $9.7122\times 10^{-3}$ & $3.0233 \times 10^{-2}$\\
$d^m_{256^2-64^2}$ & $1.0004 \times 10^{-2}$ & $1.4886\times 10^{-2}$ & $2.7392 \times 10^{-2}$\\
\hline
\multicolumn{4}{|c|}{$\text{$L_2(\Omega)$ error for $\eta$}$} \\
\hline
  & $m = t_1$ & $m = t_2$ & $m=t_3$\\
\hline
$d^m_{128^2}$      & $9.7196 \times 10^{-2}$ & $1.2153\times 10^{-1}$ & $5.3632 \times 10^{-1}$\\
$d^m_{256^2-64^2}$ & $7.7932 \times 10^{-2}$ & $1.4930\times 10^{-2}$ & $5.9167 \times 10^{-1}$\\
\hline
\end{tabular}
\end{center}
\caption{Comparisons of errors between an adaptive two-grid simulation ($256^2-64^2$) with AMR and a standard $128^2$ uniform grid simulation.}
\label{table_two_grid_2}
\end{table}

\subsection{3-D example}
\label{3d}

We mentioned previously in Remark \ref{rem1} that solving a $512^3$ 3-D simulation using a standard uniform grid requires memory of over 100 gigabytes space.
Using the two-grid solution strategy and dynamic AMR for the phase field variable, we can run a $512^3$ simulation with less than 20 gigabytes memory requirement.
This simulation is done in a 3-D domain $\Omega=(0,1)^3$.
We choose a uniform $64^3$ to be the grid where we store all the solutions and solve the adjoint equation.
The finest grid is an adaptive grid and if it was to become uniform, it would have the resolution of $512^3$.
The temporal domain is $(0,T]=(0,0.001]$, with a time step size $\tau = 5\times 10^{-5}$.
We use the same $\alpha$ and $\theta$ shown in Table \ref{table_parameters}.
More importantly, we choose an $\epsilon = 0.02$, with the finest grid and the domain, we can ensure there are roughly $10$ grid points in the interfical region in each axis direction.
It is worth noting that the interfical region of this simulation \textit{can not} be resolved accurately on any coarser grids than proposed here \cite{DDE}.

We define the initial shape to be a sphere
\begin{equation}
\label{equ_31}
\phi^{t=0} = \tanh\left( \frac{-\left[ \left(x-0.5\right)^2+\left(y-0.5\right)^2+\left(z-0.5\right)^2-0.25^2\right]}{\epsilon} \right)
\end{equation}
and the desired data to be an ellipsoid
\begin{equation}
\label{equ_32}
\phi_{obs} = \tanh\left( \frac{-\left[ \frac{\left(x-0.5\right)^2}{2}+\left(y-0.5\right)^2+\left(z-0.5\right)^2-0.25^2\right]}{\epsilon} \right).
\end{equation}

\begin{figure}[!htbp]
\centering
\includegraphics[width=15cm]{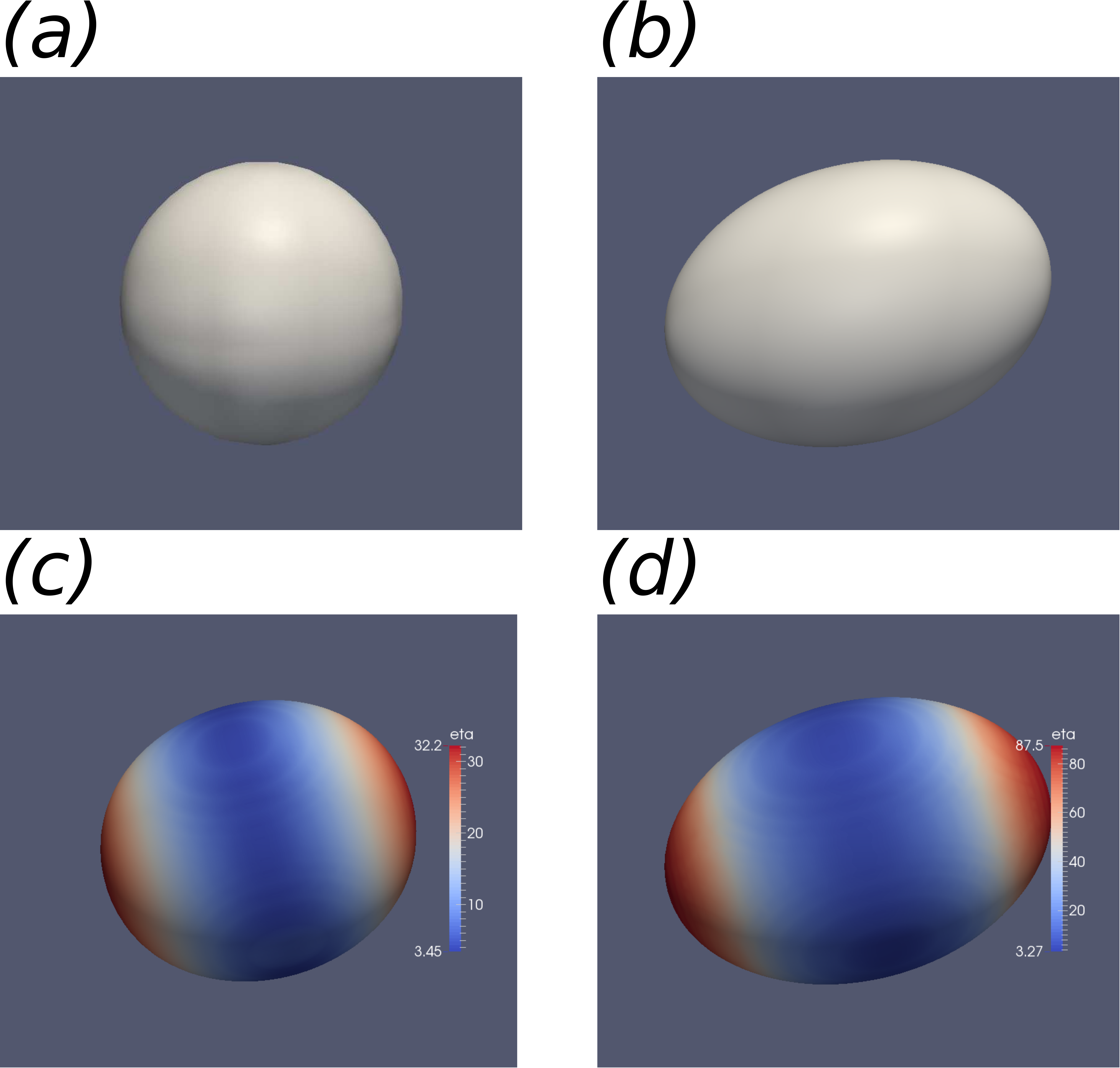}
\caption{Figures (a) and (b) show the zero-isosurface of $\phi$ of initial data (i.e. (\ref{equ_31})) and desired data (i.e. (\ref{equ_32})) respectively; (c) and (d) illustrate the zero-isosurface of computed solutions halfway through (i.e. $t=T/2$) and the final shape (i.e. $t=T=0.001$) respectively. We use colours to indicate the corresponding solutions of $\eta$ on the zero-isosurface. The colour version of this figure is online.}
\label{fig_3_D}
\end{figure}

\begin{figure}[!htbp]
\centering
\includegraphics[width=10cm]{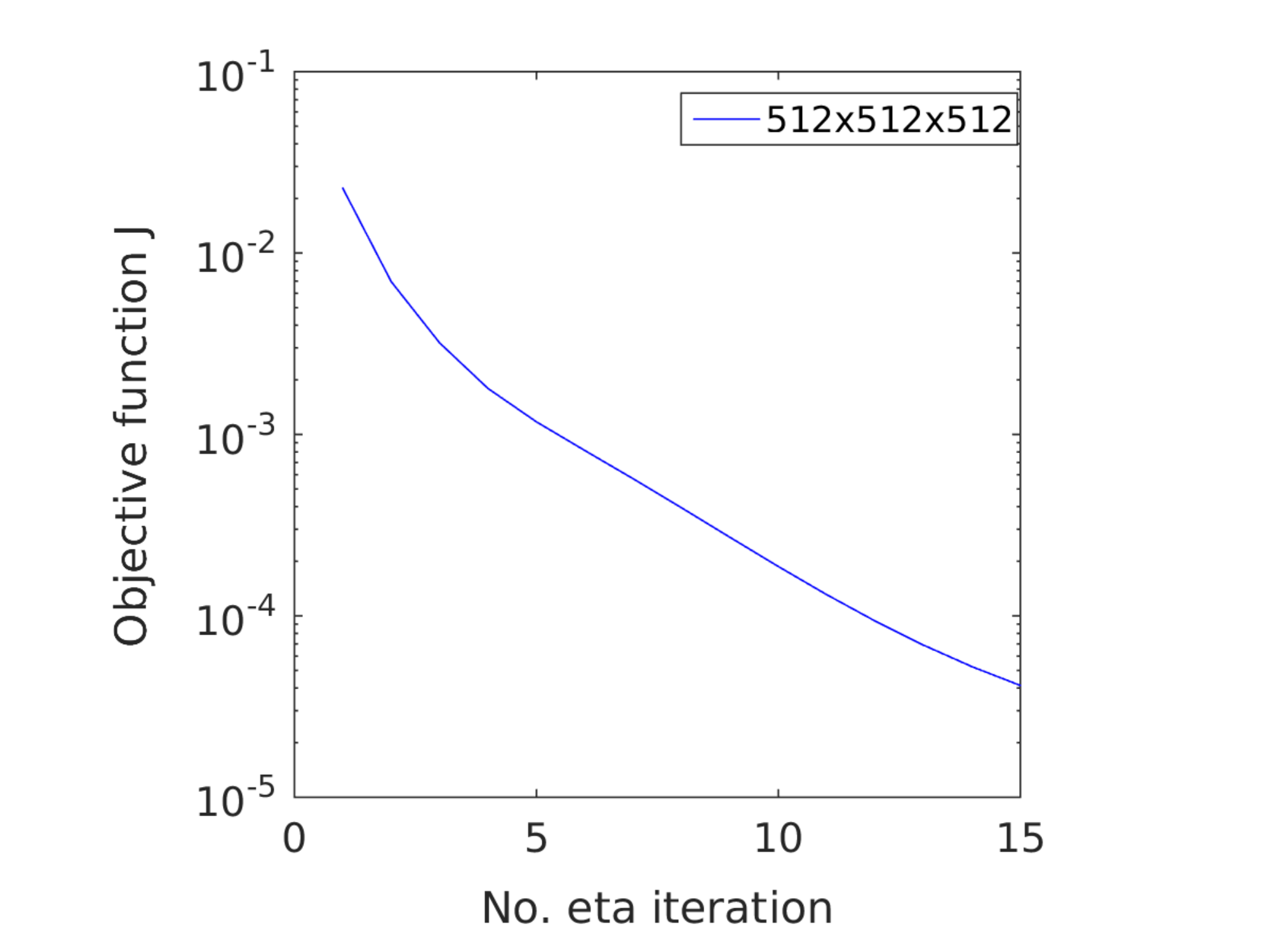}
\caption{A semi-log plot shows the reductions of the objective function $J$. The colour version of this figure is online.}
\label{fig_3_D_J}
\end{figure}

The zero-isosurface of $\phi$ for both the initial and desired shapes are illustrated in Figure \ref{fig_3_D} (a) and (b) respectively.
Following a fixed $15$ $\eta$ iterations, we present two plots of the zero-isosurface of $\phi$ in Figure \ref{fig_3_D} (c) and (d).
The solution in (c) is halfway through the temporal domain and the solution in (d) is the computed final shape.
We use colours and colour-value indicator on the side to demonstrate the corresponding solutions of $\eta$ on the zero-isosurface.
The reductions of the objective function $J$ are shown in Figure \ref{fig_3_D_J}.

\subsection{Irregular shapes}

In all our previous simulations we used relatively simple shapes for illustrative purposes only.
In this subsection, we show some irregular shapes in both 2-D and 3-D in order to illustrate that the proposed optimal control approach is capable of dealing with general interfaces.

\begin{figure}[!htbp]
\centering
\includegraphics[width=15cm]{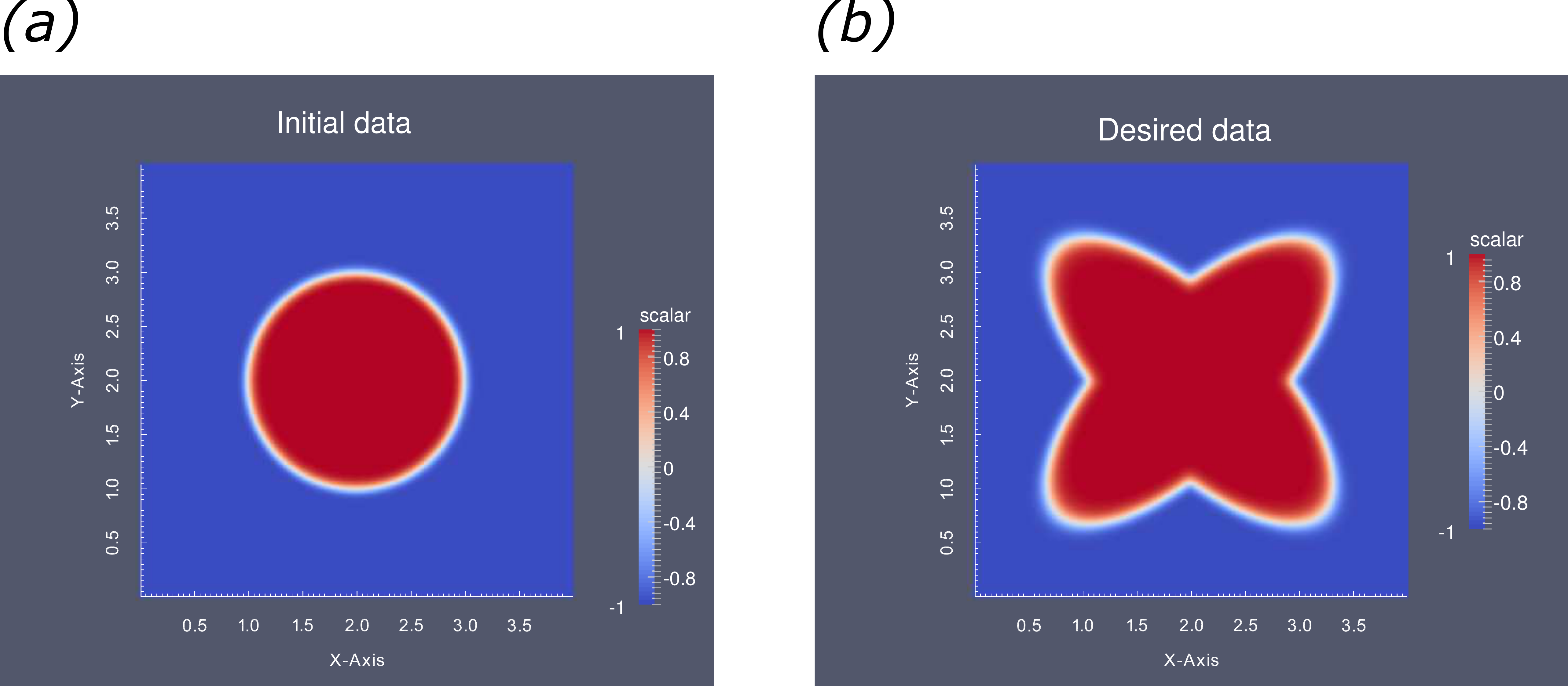}
\caption{(a) shows the initial data from \eqref{equ_2_initial} and (b) illustrates the desired data from \eqref{equ_another_desired}. The colour version of this figure is online.}
\label{fig_another_1}
\end{figure}

\begin{figure}[!htbp]
\centering
\includegraphics[width=13cm]{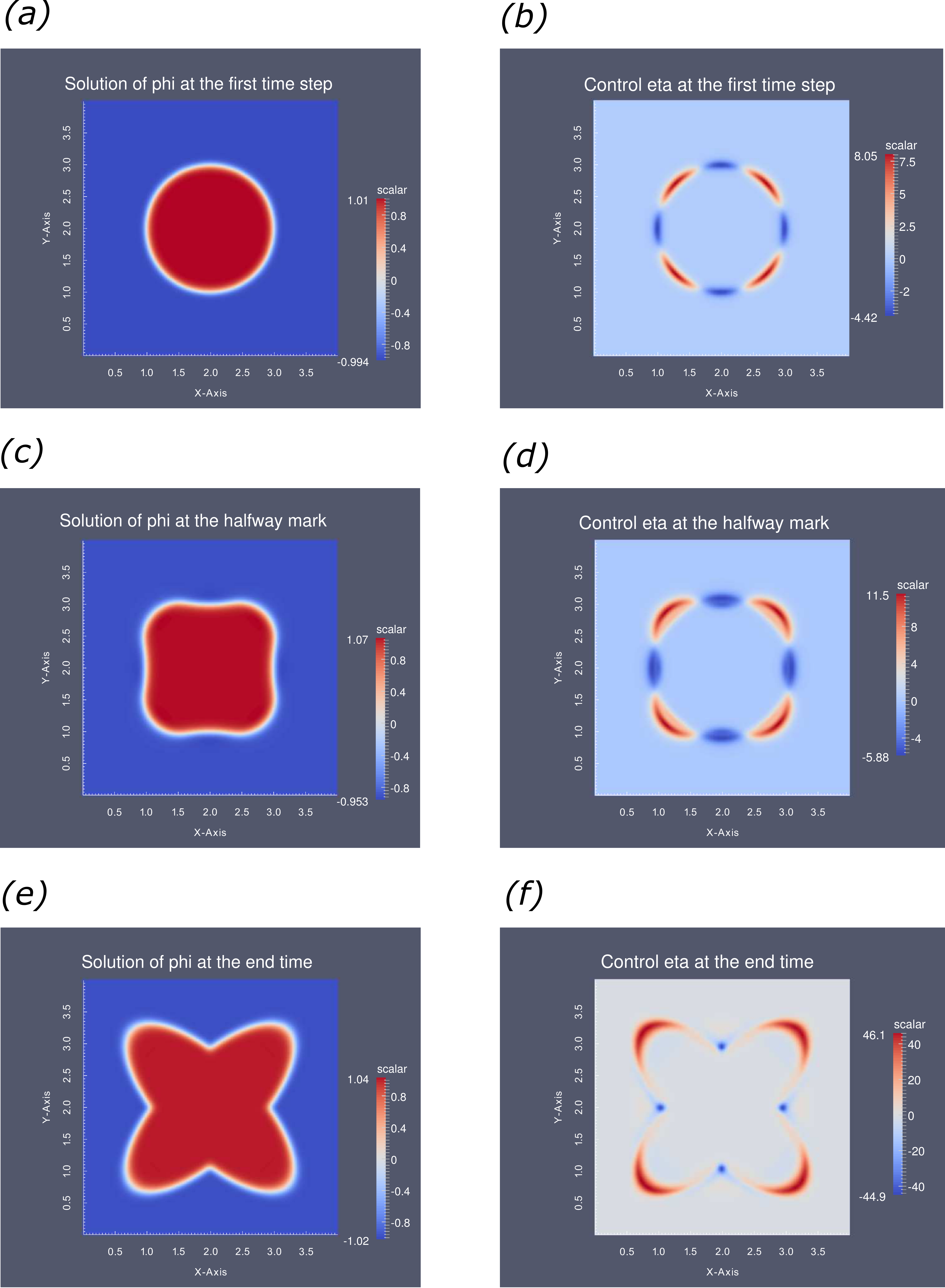}
\caption{(a) and (b) show the solutions of $\phi$ and $\eta$ at the first time step, respectively; (c) and (d) illustrate the solutions at the halfway mark (i.e. $t=0.0025$); (e) and (f) show the solutions at the final time $T$. The colour version of this figure is online.}
\label{fig_another_2}
\end{figure}

We start with a 2-D example which takes a circle as the initial shape and the desired shape is the following
\begin{equation}
\label{equ_another_desired}
\begin{split}
\phi_{obs} = \max\Bigg\{ &\tanh\left( \frac{-\left[ \frac{\left[\left(x-2\right)+\left(y-2\right)\right]^2}{6}+\frac{\left[\left(y-2\right)-\left(x-2\right)\right]^2}{1}-1\right]}{\epsilon} \right),\\
&\tanh\left( \frac{-\left[ \frac{\left[\left(x-2\right)+\left(y-2\right)\right]^2}{6}+\frac{\left[\left(y-2\right)-\left(x-2\right)\right]^2}{1}-1\right]}{\epsilon} \right) \Bigg\}.
\end{split}
\end{equation}
We take the computational domain $\Omega=(0,4)^2$ and use the parameters presented in Table \ref{table_parameters}.
This simulation is solved using a two-grid approach which has a $512^2$ grid for the Allen-Cahn and a $64^2$ for the adjoint equation.
We set $T=0.05$ and use a time step size $\tau = 0.001$.
The initial and desired data are illustrated in Figure \ref{fig_another_1}.
We present our results in Figure \ref{fig_another_2}, which include the solutions of $\phi$ at the first time step, the halfway mark (i.e. $t=0.0025$) and the end time, together with their corresponding control $\eta$.

We define two 3-D shapes as follows
\begin{equation}
\phi^0 = \tanh \left( \frac{-\left[2\left(\left(x-2\right)-\left(z-2\right)^2\right)^2 + \left(y-2\right)^2 + \left(z-2\right)^2 - 1\right]}{\epsilon} \right),
\end{equation}
\begin{equation}
\phi_{obs} = \tanh \left( \frac{-\left[\left(\left(y-2.3\right)-\left(z-2.3\right)^2\right)^2 + 2\left(x-2.3\right)^2 + \left(z-2.3\right)^2 - 1\right]}{\epsilon} \right).
\end{equation}

\begin{figure}[!htbp]
\centering
\includegraphics[width=13cm]{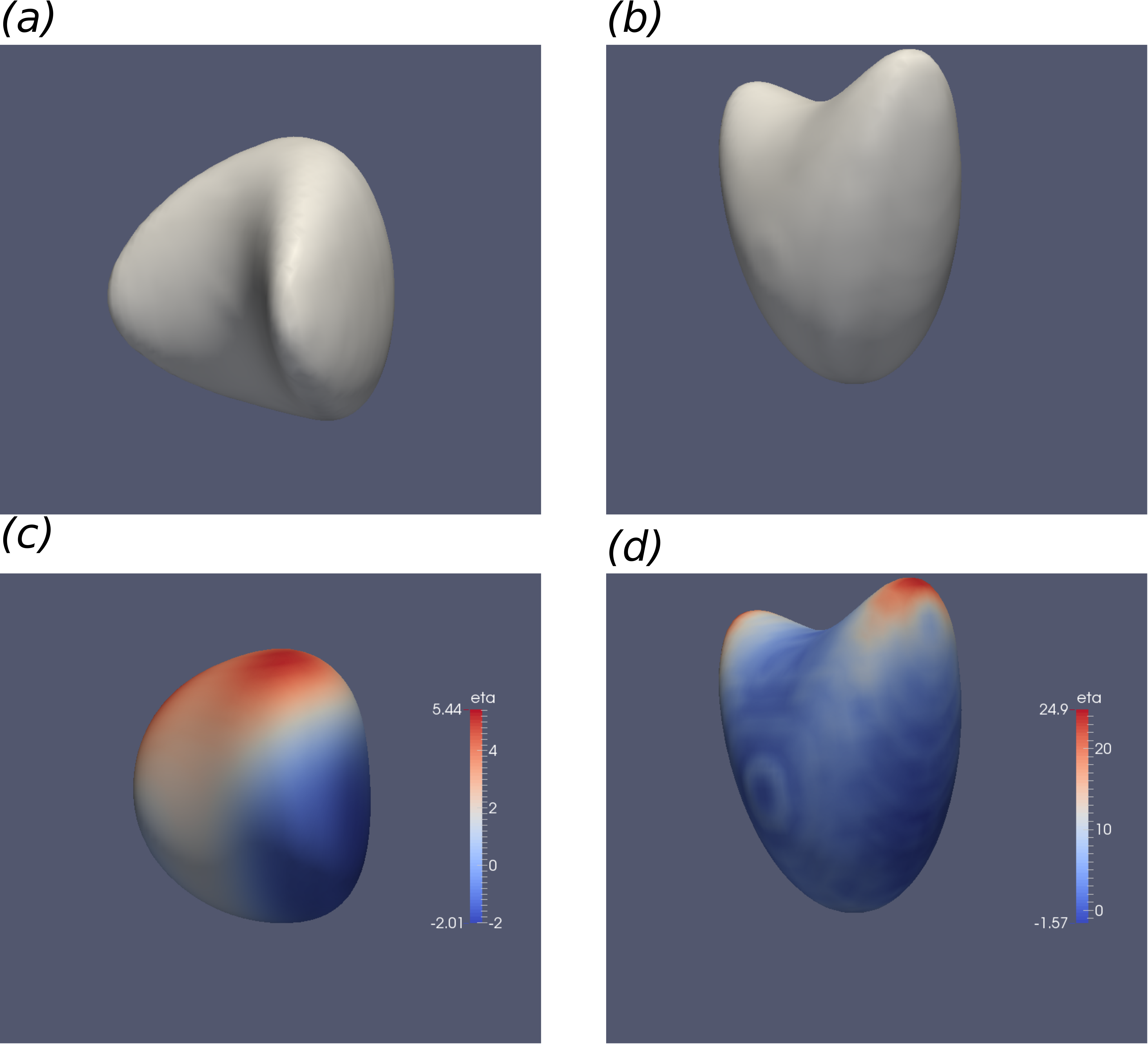}
\caption{Figures (a) and (b) show the zero-isosurface of $\phi$ of initial data and desired data respectively; (c) and (d) illustrate the zero-isosurface of computed solutions halfway through (i.e. $t=T/2$) and the final shape (i.e. $t=T=0.001$) respectively. We use colours to indicate the corresponding solutions of $\eta$ on the zero-isosurface. The colour version of this figure is online.}
\label{3d_2}
\end{figure}

The simulation has the same setting as the one described in Subsection \ref{3d} and we illustrate the zero-isosurface together with the values of the optimal control $\eta$ on this isosurface in Figure \ref{3d_2}.

\section{Conclusion}
\label{sec_concolusions}

In this work, we focussed on the development of robust and efficient solution procedures for the approximation of  the optimal control of geometric evolution laws using phase field formulations, the problems under consideration arise naturally in many applications \cite{papadakis2007variational, Hauber1,Hauber2,Troltzsch,hinze2009optimization}.
Such optimal control problems are very computationally-demanding and memory hungry especially when posed in three dimensions.
Thus the development of an efficient, robust and accurate solver is of much importance.
We have described, in detail, a solution procedure that combines a number of state-of-the-art algorithms to improve overall efficiency.
We employed a steepest descent approach for the iterative computation of the optimal control.
We introduced an adaptive-step-size algorithm which tries to use as large a step size as possible to reduce the number of iterations needed.
Robust and efficient solvers for both the forward (Allen-Cahn) and adjoint equations, based on FAS multigrid methods with MLAT are described together with their parallel implementation which is crucial for minimising wall clock time due to the massive memory requirements.
We discussed the use of mesh refinement which dramatically reduces the number of degrees of freedom required for the solution of the forward problem, and is crucial in terms of reducing the computational complexity.
A major finding of this work is that a two-grid solution strategy, in which the forward equation is solved on an adaptively refined grid whilst the adjoint problem is solved on a coarser grid, thus significantly reducing CPU and memory requirements, appears to lead to only a minor loss in accuracy.
We have implemented our algorithms and conducted detailed tests and benchmarks of our solution methods using 2-D and 3-D examples.
The conclusion is that our solution algorithms can significantly improve efficiency while maintaining an acceptable accuracy.

Possible further work, which is the focus of our current work is the application of the methodologies described in this article to real world problems. In particular, as investigated in Blazakis et al. \cite{Costas}, the optimal control problem we solve may be useful for whole cell tracking and reconstruction of dynamic cell morphologies from static imaging data. A particular advantage of the present approach is that, in contrast to the majority of existing whole cell tracking algorithms, aspects of the physics of cell migration may be encoded in the forward model and hence in the recovered trajectories. Using the techniques introduced in this paper, we may consider 3-D examples of cell tracking as well as other scenarios involving tracking multiple cells over long time horizons which requires high spatial resolution and integration over a large time interval. Our solution methods are not restricted to  forward models involving phase field formulations of geometric evolution laws; we expect that  our solution methodologies are likely to be a robust and efficient option for problems involving the control of semilinear parabolic PDEs in general. 

\subsection*{Data Management} All the computational data output is included in the present manuscript.

\section*{Acknowledgements}
All authors acknowledge support from the Leverhulme Trust Research Project Grant (RPG-2014-149). 
The work of CV, VS and AM was partially supported by the Engineering and Physical Sciences Research Council, UK grant (EP/J016780/1). 
This work (AM) has also received funding from the European Union's Horizon 2020 research and innovation programme under the Marie Sklodowska-Curie grant agreement No 642866.
The work of CV is partially supported by an EPSRC Impact Accelerator Account award. 
This research was finalised whilst all authors were participants in the Isaac Newton Institute Program, Coupling Geometric PDEs with Physics for Cell Morphology, Motility and Pattern Formation.

%
%
\bibliographystyle{unsrt}
\bibliography{refs}

\end{document}